\documentclass{article}
\usepackage{amsmath,amsthm}
\usepackage{graphicx,epstopdf,epsfig,multirow,epic,bm}
\usepackage{lineno,hyperref}
\usepackage{amsfonts}
\usepackage{amssymb}

\normalsize \rm
\begin{document}

\begin{center}
{\Large  \textbf { Mean Hitting Time on Recursive Growth Tree Network }}\\[12pt]
{\large Fei Ma$^{a,}$\footnote{~The author's E-mail: mafei123987@163.com. },\quad  Ping Wang$^{b,c,d,}$\footnote{~The author's E-mail: pwang@pku.edu.cn.} }\\[6pt]
{\footnotesize $^{a}$ School of Electronics Engineering and Computer Science, Peking University, Beijing 100871, China\\
$^{b}$ National Engineering Research Center for Software Engineering, Peking University, Beijing, China\\
$^{c}$ School of Software and Microelectronics, Peking University, Beijing  102600, China\\
$^{d}$ Key Laboratory of High Confidence Software Technologies (PKU), Ministry of Education, Beijing, China}\\[12pt]
\end{center}

\begin{quote}
\textbf{Abstract:} In this paper, we are concerned with mean hitting time $\langle\mathcal{H}\rangle$ for random walks on recursive growth tree networks that are built based on an arbitrary tree as the seed via implementing various primitive graphic operations, and propose a series of combinatorial techniques that are called Mapping Transformation to exactly determine the associated $\langle\mathcal{H}\rangle-$polynomial. Our formulas can be able to completely cover the previously published results in some well-studied and specific cases where a single edge or a star is often chose to serve as seed for creating recursive growth models. The techniques proposed are more convenient than the commonly-used spectral methods mainly because of getting around the operations of matrix inversion and multiplication. Accordingly, our results can be extended for both many other stochastic models including BA-scale-free tree and random uniform tree as well as graphs of great interest consisting of line graph of tree and Vicsek fractal network to derive numerical solutions of related structural parameters. And then, the closed-form solutions of two extensions of Wiener index with respect to multiplicative and additive degrees on an arbitrary tree are conveniently obtained as well. In addition, we discuss some extremal problems of random walks on tree networks and outline the related research directions in the next step.
\\


\end{quote}

\section{Introduction}

As is known to us all, various kinds of dynamics taking place on many complex systems in nature and society can be well described as random walks on complex network, an abstract mathematical object that has attracted considerable attention in science fields including applied mathematics, theoretical computer sciences, and statistical physics in the past \cite{L-L-1993}-\cite{Sylvester-2021}. Studied examples include stochastic behaviour of molecules in rarified gases \cite{Dongari-2011}, protein folding and misfolding \cite{Dobson-2003}, information flow in social networks \cite{Guille-2012}, traffic and mobility patterns on the internet \cite{Gonzalez-2008}, fluctuations in stock prices \cite{Wang-2010}, information dissemination on space \cite{Lam-2012} as well as the behaviour of stochastic search algorithms \cite{Dorigo-2005}, and so forth. In the context of random walks on networks, one of most important and fundamental structural parameters is the so-called hitting time \cite{Oliver-2013}. It is the expected time taken by a walker performing random walks on network to first reach to its destination vertex. Indeed, it is using such a parameter that enables ones to well understand this paradigmatic dynamical process and to uncover the effects of the underlying topological structure of network on the behavior of random-walk dynamics \cite{Volchenkov-2011}-\cite{Guerin-2016}.

In general, hitting time from source vertex to destination vertex is analytically estimated based on the spectral of Laplacian matrix corresponding to network \cite{Zhu-1996}. Accordingly, mean hitting time for a network in question, which is defined as the average over hitting times of all possible pairs of vertices on network, is exactly obtained as shown in Lemma 2.1. In theory, this is a versatile tool suitable for an arbitrary network, and enables exact estimation of mean hitting time. On the other hand, such a calculation technique is built based on some matrix manipulations, such as matrix inversion and multiplication, which often has a higher complexity. For example, exactly computing the invert of a matrix needs to take $O(n^{3})$ time. This is computationally impractical for large network with $n$ vertices. In recently published paper \cite{Peng-2021}, Peng \emph{et al} have substantially improved algorithm for performing matrix inversion from time complexity $O(n^{3})$ to $O(n^{2.3316})$. Indeed, it is a great progress for obtaining appropriate solutions of related parameters based on fundamental matrix arithmetics from the theory point of view. Nonetheless, perhaps it is still not an available manner in which one wants to obtain the closed-from solutions of those structural parameters, such as mean hitting time, on some large networks. Therefore, the need increases for a specific and mathematically rigorous tool that addresses this kind of problems on networks of great interest, such as tree network. Based on this, the goal of this paper is to investigate random walks on various tree networks, which have a wide range of applications in many areas, such as date center networks \cite{Couto-2016} and material science \cite{Blumen-2004}, and we establish more effective methods to derive exact solutions of relevant parameters including mean hitting time. Additionally, many other interesting graphs are also chosen as our focuses. Note that, throughout this paper, the terms graph and network are used indistinctly for convenience.

Tree, denoted by $\mathcal{T}$, as the simplest and most fundamental connected graph model, has been widely studied in the rich literature, especially, in graph theory and theoretical computer science \cite{Cooper-2010}-\cite{Georgakopoulos-2017}, and is commonly used to model the underlying structure of some dendrimers and regular hyperbranched polymers in both bio-chemistry and material science \cite{Blumen-2004,Gurtovenko-2005}. In addition, it is well known that tree is a frequently used architecture in computer science, such as $R-$branch tree \cite{Arge-2008}. As a result, random walks on tree have also been discussed in detail \cite{Cooper-2010}-\cite{Georgakopoulos-2017} for the purpose of probing some intriguing topological properties. For example, Baronchelli \emph{et al} in \cite{Baronchelli-2008} derived the scaling of some related quantities on complex tree networks, and concluded that the absence of loops reflects in physical observables showing large differences with respect to their looped counterparts. In the prior work \cite{Ma-pre-2020}, a close relationship between mean hitting time and Wiener index on tree is established using fundamental matrix operation. Again for instance, Beveridge \emph{et al} in \cite{Beveridge-2013} characterize the extremal structures for certain random walks on trees and obtain some inequalities satisfied by related parameters including mean hitting time. In \cite{Georgakopoulos-2017}, Georgakopoulos \emph{et al} exhibited a close connection between hitting times of the simple random walk on a tree, Wiener index, and related graph invariants.

Clearly, most prior work focuses particularly on revealing some general principles associated with structural parameters for random walks on tree. Besides that, some specific tree networks of great interest, such as T-graph \cite{Redner-2001}, Cayley trees \cite{Wu-2012}, and Vicsek fractals as models of polymer networks \cite{Vicsek-1983}, have also received more attention for special-purpose. Accordingly, the related parameters including mean hitting time have been in depth studied, and the corresponding precise solutions were obtained. It should be mentioned that all almost solutions are derived based on the commonly used spectral technique (i.e., Lemma 2.1). At the same time, all the specific tree models mentioned above are recursively built up given a single edge or a star as the seed. Indeed, the previously published results in \cite{Cooper-2010}-\cite{Georgakopoulos-2017} have provided ones with good understanding on structure of the resultant tree networks in the limit of large graph size. On the other hand, if we are not limited to the above assumption, then it is of interest to derive more general formula of mean hitting time from the mathematical point of view. In other words, so far, the problem of determining closed-form expression of mean hitting time on recursive growth tree networks whose seed is an arbitrary tree is still open.

Therefore, in this work, we aim to completely address the issue above. Towards this end, we come up with a series of more effective combinatorial methods, which are called Mapping Transformation hereafter. The reason for this is two-fold. (i) As will see in the following, the development of these methods is closely dependent on the nature of tree itself. This highly reflects the principle of Ockham's razor. (ii) It is nontrivial to extend some previous methods, which are mainly based upon Laplacian matrix, suitable for the simplest case to a more general situation as considered here. For example, it has been shown in \cite{Jayanthi-1992} that even for the simplest case of Vicsek fractal, some complicated and advanced mathematical tools such as real space Green's function approach have to be adopted to derive closed-form solution of mean hitting time when making use of spectral technique. Intuitively, facing with more general situation where an arbitrary tree is chose as the seed as discussed in the rest of this paper, this kind of popularly-used techniques can become prohibitively complicated. In a nutshell, we propose more manageable methods in order to overcome the above challenge. Note also that the proposed methods are also applicable for many other stochastic graphs to measure some structural parameters.

Before beginning with detailed demonstrations, we first review the main results, which are organized into the following theorems and propositions. See Sections 2-5 for a more thorough treatment of the necessary background. Without loss of generality, we assume that seed $\mathcal{T}$ has $n$ vertices and its Wiener index is $\mathcal{W}_{\mathcal{T}}$ in the remainder of this paper. In view of pattern-specific of parameter $m$, the correspondingly concrete definition is deferred to show in Section 3.

\textbf{Theorem I.1} \emph{The exact solution of Winer index $\mathcal{W}_{\mathcal{T}^{m}(1)}$ on tree $\mathcal{T}^{m}(1)$ is given by}

\begin{equation}\label{eqa:MF-4-1-0}
\mathcal{W}_{\mathcal{T}^{m}(1)}=(m+1)^{3}\mathcal{W}_{\mathcal{T}}-(n-1)(m+1)m^{2}-2\left[m\left(
                                                             \begin{array}{c}
                                                               n-1 \\
                                                               2 \\
                                                             \end{array}
                                                           \right)+\left(
                                                             \begin{array}{c}
                                                               n \\
                                                               2 \\
                                                             \end{array}
                                                           \right)\right]\left(
                                                             \begin{array}{c}
                                                               m+1 \\
                                                               2 \\
                                                             \end{array}
                                                           \right)+(n-1)\left(
                                                             \begin{array}{c}
                                                               m+1 \\
                                                               3 \\
                                                             \end{array}
                                                           \right).
\end{equation}

\textbf{Theorem I.2} \emph{The exact solution of Winer index $\mathcal{W}_{\mathcal{T}^{\mathrm{I}}_{m}(1)}$ on tree $\mathcal{T}^{\mathrm{I}}_{m}(1)$ is given by}

\begin{equation}\label{eqa:MF-4-2-0}
\mathcal{W}_{\mathcal{T}^{\mathrm{I}}_{m}(1)}=(m+1)^{2}\mathcal{W}_{\mathcal{T}}+2m(m+1)\left(
                                                             \begin{array}{c}
                                                               n \\
                                                               2 \\
                                                             \end{array}
                                                           \right)+n\left[m+2\left(
                                                             \begin{array}{c}
                                                               m \\
                                                               2 \\
                                                             \end{array}
                                                           \right)\right].
\end{equation}

\textbf{Theorem I.3} \emph{The exact solution of Winer index $\mathcal{W}_{\mathcal{T}^{\mathrm{T}}_{m}(1)}$ on tree $\mathcal{T}^{\mathrm{T}}_{m}(1)$ is given by}

\begin{equation}\label{eqa:MF-4-3-0}
\mathcal{W}_{\mathcal{T}^{\mathrm{T}}_{m}(1)}=2(m+2)^{2}\mathcal{W}_{\mathcal{T}}-(n-1)\left[2m^{2}+3m+2n-2\left(
                                                             \begin{array}{c}
                                                               m \\
                                                               2 \\
                                                             \end{array}
                                                           \right)\right]-2m\left(
                                                             \begin{array}{c}
                                                               n-1 \\
                                                               2 \\
                                                             \end{array}
                                                           \right).
\end{equation}

\textbf{Theorem I.4} \emph{The exact solution of Winer index $\mathcal{W}_{\mathcal{T}^{\mathrm{V}}_{m}(1)}$ on tree $\mathcal{T}^{\mathrm{V}}_{m}(1)$ is given by}

\begin{equation}\label{eqa:MF-4-4-0}
\mathcal{W}_{\mathcal{T}^{\mathrm{V}}_{m}(1)}=3(m+1)^{2}\mathcal{W}_{\mathcal{T}}+nm^{2}+2[(m-1)^{2}+m-3]\left(
                                                             \begin{array}{c}
                                                               n \\
                                                               2 \\
                                                             \end{array}
                                                           \right).
\end{equation}

\textbf{Theorem I.5} \emph{The exact solution of Winer index $\mathcal{W}_{\mathcal{T}^{\mathrm{II}}_{m}(1)}$ on tree $\mathcal{T}^{\mathrm{II}}_{m}(1)$ is given by}

\begin{equation}\label{eqa:MF-4-5-0}
\begin{aligned}\mathcal{W}_{\mathcal{T}^{\mathrm{II}}_{m}(1)}&=(2m+1)^{2}\mathcal{W}_{\mathcal{T}}-(n-1)(2n-1)m^{2}+\frac{8m^{2}}{3}(n-1)^{2}-2m(n-1)\\
&\quad+\frac{4m}{3}n(n-1)-\frac{2m}{3}\left(
                                                             \begin{array}{c}
                                                               n \\
                                                               2 \\
                                                             \end{array}
                                                           \right)+\frac{4m^{2}}{3}\left[\left(
                                                             \begin{array}{c}
                                                               n \\
                                                               2 \\
                                                             \end{array}
                                                           \right)+\left(
                                                             \begin{array}{c}
                                                               n-1 \\
                                                               2 \\
                                                             \end{array}
                                                           \right)\right]
                                                           \end{aligned}.
\end{equation}

\textbf{Theorem I.6} \emph{The exact solution of Winer index $\mathcal{W}_{\mathcal{T}^{\mathrm{III}}_{m}(1)}$ on tree $\mathcal{T}^{\mathrm{III}}_{m}(1)$ is given by}

\begin{equation}\label{eqa:MF-4-6-0}
\begin{aligned}\mathcal{W}_{\mathcal{T}^{\mathrm{III}}_{m}(1)}&=(m-1)^{2}\mathcal{W}_{\mathcal{T}}+n(n-1)\left(m^{2}-\frac{m}{3}-\frac{4}{3}\right)+nm^{2}-2(n-1)(m-1)\\
&\quad-(n-1)(2n-1)+\left(\frac{8}{3}-2m\right)(n-1)^{2}+\frac{2(m+1)}{3}\left(
                                                             \begin{array}{c}
                                                               n \\
                                                               2 \\
                                                             \end{array}
                                                           \right)+\frac{4}{3}\left[\left(
                                                             \begin{array}{c}
                                                               n \\
                                                               2 \\
                                                             \end{array}
                                                           \right)+\left(
                                                             \begin{array}{c}
                                                               n-1 \\
                                                               2 \\
                                                             \end{array}
                                                           \right)\right]
\end{aligned},
\end{equation}

Using Theorems I.1-I.6, we can show the potential applications, i.e., deriving the numerical solutions of relevant structural parameters on many widely-studied network models including BA-scale-free tree \cite{Albert-1999-1}, line graph of tree \cite{Knor-2013}, random uniform tree \cite{Chelminiak-2011}, and Vicsek fractal network \cite{Vicsek-1983}, which are as follows.

\textbf{Proposition I.7} \emph{The closed-form solution of mean hitting time $\langle\mathcal{H}_{\mathcal{T}^{m}(t)}\rangle$ of subdivision tree $\mathcal{T}^{m}(t)$ is given by}

\begin{equation}\label{eqa:MF-5-1-0}
\langle\mathcal{H}_{\mathcal{T}^{m}(t)}\rangle=\frac{2}{|\mathcal{T}^{m}(t)|}\left\{\begin{aligned}&[f_{1}(m)]^{t}\mathcal{W}_{\mathcal{T}}-g_{1}(m)\sum_{i=0}^{t-1}[f_{1}(m)]^{i}|\mathcal{T}^{m}(t-1-i)|^{2}\\
&+h_{1}(m)\sum_{i=0}^{t-1}[f_{1}(m)]^{i}|\mathcal{T}^{m}(t-1-i)|-l_{1}(m)\sum_{i=0}^{t-1}[f_{1}(m)]^{i}
\end{aligned}\right\},
\end{equation}
\emph{in which $f_{1}(m)=(m+1)^{3}$, $g_{1}(m)=\frac{m(m+1)^{2}}{2}$, $h_{1}(m)=\frac{m(m+1)(2m+1)}{3}$, $l_{1}(m)=\frac{m(m^{2}-1)}{6}$ and $|\mathcal{T}^{m}(t)|=(n-1)(m+1)^{t}+1$.}

\textbf{Proposition I.8} \emph{The closed-form solution of Wiener index $\mathcal{W}_{\mathcal{L}}$ of line graph $\mathcal{T}_{\mathcal{L}}$ of a given tree $\mathcal{T}$ is given by}
\begin{equation}\label{eqa:MF-5-2-0}
\mathcal{W}_{\mathcal{L}}=\mathcal{W}_{\mathcal{T}}-\left(
                                                             \begin{array}{c}
                                                               n \\
                                                               2 \\
                                                             \end{array}
                                                           \right).
\end{equation}

\textbf{Proposition I.9} \emph{The analytic solution of mean shortest path length $\langle\mathcal{W}_{SF}(t)\rangle$ of the classic BA-scale-free tree $\mathcal{T}_{SF}(t)$ is given by}

\begin{equation}\label{eqa:MF-5-3-0}
\langle\mathcal{W}_{SF}(t)\rangle=\frac{2}{(t+2)(t+1)}\left[\frac{t+1}{2}-\frac{2t+1}{4t^{2}}+4\prod_{i=2}^{t}\left(\frac{i+1}{i}\right)^{2}+\sum_{i=2}^{t-1}\left(\frac{i+1}{2}-\frac{2i+1}{4i^{2}}\right)\prod_{j=i+1}^{t}\left(\frac{j+1}{j}\right)^{2}\right].
\end{equation}

\textbf{Proposition I.10} \emph{The closed-form solution of mean hitting time $\langle\mathcal{H}_{\mathcal{T}^{\mathrm{II}}_{m}(t)}\rangle$ of scale-free trees $\mathcal{T}^{\mathrm{II}}_{m}(t)$ is given by}

\begin{equation}\label{eqa:MF-5-4-0}
\langle\mathcal{H}_{\mathcal{T}^{\mathrm{II}}_{m}(t)}\rangle=\frac{2}{|\mathcal{T}^{\mathrm{II}}_{m}(t)|}\left\{\begin{aligned}&[f_{2}(m)]^{t}\mathcal{W}_{\mathcal{T}}+g_{2}(m)\sum_{i=0}^{t-1}[f_{2}(m)]^{i}|\mathcal{T}^{\mathrm{II}}_{m}(t-1-i)|^{2}\\
&-h_{2}(m)\sum_{i=0}^{t-1}[f_{2}(m)]^{i}|\mathcal{T}^{\mathrm{II}}_{m}(t-1-i)|+l_{2}(m)\sum_{i=0}^{t-1}[f_{2}(m)]^{i}
\end{aligned}\right\},
\end{equation}
\emph{in which $f_{2}(m)=(2m+1)^{2}$, $g_{2}(m)=m(2m+1)$, $h_{2}(m)=m(5m+3)$, $l_{2}(m)=m(3m+2)$ and $|\mathcal{T}^{\mathrm{II}}_{m}(t)|=(n-1)(2m+1)^{t}+1$.}

\textbf{Proposition I.11} \emph{The analytic solution of mean shortest path length $\langle\mathcal{W}_{\mathcal{T}_{RG}}(t)\rangle$ of random uniform growth trees $\mathcal{T}_{RG}(t)$ is given by}

\begin{equation}\label{eqa:MF-5-5-0}
\langle\mathcal{W}_{\mathcal{T}_{RG}}(t)\rangle=\frac{2}{(t+2)(t+1)}\left[t+1+\frac{t}{t+1}+\prod_{i=0}^{t}\left(\frac{i+2}{i+1}\right)+\sum_{i=1}^{t-1}\left(i+1+\frac{i}{i+1}\right)\prod_{j=i+1}^{t}\left(\frac{i+2}{i+1}\right)\right].
\end{equation}

\textbf{Proposition I.12} \emph{The closed-form solution of mean hitting time $\langle\mathcal{H}_{\mathcal{T}^{\mathrm{I}}_{m}(t)}\rangle$ of exponential tree $\mathcal{T}^{\mathrm{I}}_{m}(t)$ is given by}

\begin{equation}\label{eqa:MF-5-6-0}
\langle\mathcal{H}_{\mathcal{T}^{\mathrm{I}}_{m}(t)}\rangle=\frac{2}{|\mathcal{T}^{\mathrm{I}}_{m}(t)|}\left\{\begin{aligned}&[f_{3}(m)]^{t}\mathcal{W}_{\mathcal{T}}+g_{3}(m)\sum_{i=0}^{t-1}[f_{3}(m)]^{i}|\mathcal{T}^{\mathrm{I}}_{m}(t-i-1)|^{2}\\
&-h_{3}(m)\sum_{i=0}^{t-1}[f_{3}(m)]^{i}|\mathcal{T}^{\mathrm{I}}_{m}(t-i-1)|
\end{aligned}\right\},
\end{equation}
\emph{in which $f_{3}(m)=(m+1)^{2}$, $g_{3}(m)=m(m+1)$, $h_{3}(m)=m$ and $|\mathcal{T}^{\mathrm{I}}_{m}(t)|=n(m+1)^{t}$.}

\textbf{Proposition I.13} \emph{The closed-form solution of mean hitting time $\langle\mathcal{H}_{\mathcal{T}^{\mathrm{T}}_{m}(t)}\rangle$ of generalized T-fractal $\mathcal{T}^{\mathrm{T}}_{m}(t)$ is given by}

\begin{equation}\label{eqa:MF-5-7-0}
\langle\mathcal{H}_{\mathcal{T}^{\mathrm{T}}_{m}(t)}\rangle=\frac{2}{|\mathcal{T}^{\mathrm{T}}_{m}(t)|}\left\{\begin{aligned}&[f_{4}(m)]^{t}\mathcal{W}_{\mathcal{T}}-g_{4}(m)\sum_{i=0}^{t-1}[f_{4}(m)]^{i}|\mathcal{T}^{\mathrm{T}}_{m}(t-1-i)|^{2}\\
&-h_{4}(m)\sum_{i=0}^{t-1}[f_{4}(m)]^{i}|\mathcal{T}^{\mathrm{T}}_{m}(t-1-i)|+l_{4}(m)\sum_{i=0}^{t-1}[f_{4}(m)]^{i}
\end{aligned}\right\},
\end{equation}
\emph{in which $f_{4}(m)=2(m+2)^{2}$, $g_{4}(m)=(m+2)$, $h_{4}(m)=(m-1)(m+2)$, $l_{4}(m)=m^{2}+2m$ and $|\mathcal{T}^{\mathrm{T}}_{m}(t)|=(n-1)(m+2)^{t}+1$.}

\textbf{Proposition I.14} \emph{The closed-form solution of mean hitting time $\langle\mathcal{H}_{\mathcal{T}^{\mathrm{V}}_{m}(t)}\rangle$ of generalized V-fractal $\mathcal{T}^{\mathrm{V}}_{m}(t)$ is given by}

\begin{equation}\label{eqa:MF-5-8-0}
\langle\mathcal{H}_{\mathcal{T}^{\mathrm{V}}_{m}(t)}\rangle=\frac{2}{|\mathcal{T}^{\mathrm{V}}_{m}(t)|}\left\{\begin{aligned}&[f_{5}(m)]^{t}\mathcal{W}_{\mathcal{T}}+g_{5}(m)\sum_{i=0}^{t-1}[f_{5}(m)]^{i}|\mathcal{T}^{\mathrm{V}}_{m}(t-1-i)|^{2}\\
&+h_{5}(m)\sum_{i=0}^{t-1}[f_{5}(m)]^{i}|\mathcal{T}^{\mathrm{V}}_{m}(t-1-i)|
\end{aligned}\right\},
\end{equation}
\emph{in which $f_{5}(m)=3(m+1)^{2}$, $g_{5}(m)=(m-1)(m+2)$, $h_{5}(m)=(m+2)$ and $|\mathcal{T}^{\mathrm{V}}_{m}(t)|=n(m+1)^{t}$.}

\textbf{Proposition I.15} \emph{The closed-form solution of mean hitting time $\langle\mathcal{H}_{\mathcal{T}^{\mathrm{III}}_{m}(t)}\rangle$ of generalized Cayley tree $\mathcal{T}^{\mathrm{III}}_{m}(t)$ is given by}

\begin{equation}\label{eqa:MF-5-9-0}
\langle\mathcal{H}_{\mathcal{T}^{\mathrm{III}}_{m}(t)}\rangle=\frac{2}{|\mathcal{T}^{\mathrm{III}}_{m}(t)|}\left\{\begin{aligned}&[f_{6}(m)]^{t}\mathcal{W}_{\mathcal{T}}+g_{6}(m)\sum_{i=0}^{t-1}[f_{6}(m)]^{i}|\mathcal{T}^{\mathrm{III}}_{m}(t-1-i)|^{2}\\
&+h_{6}(m)\sum_{i=0}^{t-1}[f_{6}(m)]^{i}|\mathcal{T}^{\mathrm{III}}_{m}(t-1-i)|+l_{6}(m)\sum_{i=0}^{t-1}[f_{6}(m)]^{i}
\end{aligned}\right\},
\end{equation}
\emph{in which $f_{6}(m)=(m-1)^{2}$, $g_{6}(m)=(m-1)^{2}$, $h_{6}(m)=2(m-1)$, $l_{6}(m)=1$ and $|\mathcal{T}^{\mathrm{III}}_{m}(t)|=(n+\frac{2}{m-2})(m-1)^{t}-\frac{2}{m-2}$.}

\textbf{Proposition I.16} \emph{The closed-form solution of multiplicative degree Wiener index $\mathcal{W}^{\ast}_{\mathcal{T}}$ of an arbitrary tree $\mathcal{T}$ is given by}

\begin{equation}\label{eqa:MF-5-10-0}
\mathcal{W}^{\ast}_{\mathcal{T}}=\frac{1}{2}\sum_{u,v}k_{u}k_{v}d_{uv}=4\mathcal{W}_{\mathcal{T}}-(n-1)(2n-1),
\end{equation}
\emph{where $\mathcal{W}_{\mathcal{T}}$ is Wiener index of tree $\mathcal{T}$.}

\textbf{Proposition I.17} \emph{The closed-form solution of additive degree Wiener index $\mathcal{W}^{\dagger}_{\mathcal{T}}$ of an arbitrary tree $\mathcal{T}$ is given by}

\begin{equation}\label{eqa:MF-5-11-0}
\mathcal{W}^{\dagger}_{\mathcal{T}}=\frac{1}{2}\sum_{u,v}(k_{u}+k_{v})d_{uv}=4\mathcal{W}_{\mathcal{T}}-n(n-1),
\end{equation}
\emph{where $\mathcal{W}_{\mathcal{T}}$ is Wiener index of tree $\mathcal{T}$.}

The rest of this paper is organized into the next several sections. Section 2 introduces some basic terminologies, such as graph and its matrix representation, Wiener index, and random walks on graph. Section 3 presents several different types of graphic operations that will be used to create recursive growth tree networks. Section 4 shows rigorous proofs of our main results listed out in Theorems I.1-I.6, i.e., analytically determining the exact solutions of Wiener index on growth trees built in the preceding section by developing a series of combinatorial manners called Mapping Transformation. Section 5 elaborates on various applications of the results derived in Section 4 to some well-known tree networks, namely, providing strict proofs of numerical formulas for related structural parameters shown in Proposition I.9-I.17. Section 6 discusses some extremal problems and future research directions. Finally, we close this paper in Section 7.

\section{Terminologies}

In this section, we will introduce some basic concepts and notations for graphs and random walks on graphs. For convenience, we denote by $[a,b]$ a set of integers $\{a,a+1,\dots,b\}$.

\subsection{Graph and its matrix representation \cite{Biggs-1993}}

A graph (or network) $\mathcal{G}(\mathcal{V},\mathcal{E})$ is an ordered pair ($\mathcal{V}(\mathcal{G}),\mathcal{E}(\mathcal{G})$) consisting of a set $\mathcal{V}(\mathcal{G})$ of vertices and a set $\mathcal{E}(\mathcal{G})$ of edges running between vertices. Unless otherwise specified, let $\mathcal{G}$ denote a graph for brevity. The total number of vertices is denoted by $|\mathcal{V}|$ and $|\mathcal{E}|$ represents the edge number. Hereafter, all the discussed graphs are simple, unweighted and connected, namely, without multi-edges and loops.

More generally, it is convention to interpret a graph $\mathcal{G}(\mathcal{V},\mathcal{E})$ using its adjacency matrix $\mathbf{A}_{\mathcal{G}}=(a_{ij})$ in the following form

$$a_{ij}=\left\{\begin{aligned}&1, \quad\text{vertex $i$ is adjacent to $j$}\\
&0,\quad\text{otherwise}.
\end{aligned}\right.
$$
This thus encompasses some basic information about a graph itself, such as, the degree $k_{i}$ of vertex $i$ is equal to $k_{i}=\sum_{j=1}^{|\mathcal{V}|}a_{ij}$. Also, the diagonal matrix, denoted by $\mathbf{D}_{\mathcal{G}}$, may be immediately defined as follows: the $i$th diagonal entry is $k_{i}$, while all non-diagonal entries are zero, i.e., $\mathbf{D}_{\mathcal{G}}=\text{diag}[k_{1},k_{2},\dots,k_{|\mathcal{V}|}]$.

\subsection{Wiener index on graph \cite{Bondy-2008}}

In the language of graph theory, distance of a pair of vertices $u$ and $v$, denoted by $d_{uv}$, of graph $\mathcal{G}(\mathcal{V},\mathcal{E})$ is the length of a shortest path between vertices $u$ and $v$. In some published papers \cite{Katzav-2018}, this index is also called shortest path length for this pair of vertices. For a given graph $\mathcal{G}(\mathcal{V},\mathcal{E})$ as a whole, the summation over distances $d_{uv}$ of all possible pairs of vertices $u$ and $v$ is defined as Wiener index, denoted by $\mathcal{W}_{\mathcal{G}}$, namely,

\begin{equation}\label{eqa:MF-2-2-1}
\mathcal{W}_{\mathcal{G}}=\frac{1}{2}\sum_{u,v\in\mathcal{V}}d_{uv}=\sum_{1=i<j=|\mathcal{V}|}d_{ij}.
\end{equation}
here $i$ is a unique label for each vertex. Accordingly, mean shortest path length of graph $\mathcal{G}(\mathcal{V},\mathcal{E})$ is defined as follows

\begin{equation}\label{eqa:MF-2-2-2}
\langle\mathcal{W}_{\mathcal{G}}\rangle=2\mathcal{W}_{\mathcal{G}}/|\mathcal{V}|(|\mathcal{V}|-1).
\end{equation}

\subsection{Random walks \cite{Oliver-2013}}

We now consider an unbiased discrete-time random walk taking place on graph $\mathcal{G}(\mathcal{V},\mathcal{E})$. Particularly, a walker starting out from its current location $u$ moves with a uniform probability proportional to its degree $k_{u}$ to each vertex $v$ of its neighboring set in one step \cite{Aldous-1999}. In general, such a dynamical process can be certainly represented by the transition matrix $\mathbf{P}_{\mathcal{G}}=\mathbf{D}_{\mathcal{G}}^{-1}\mathbf{A}_{\mathcal{G}}$ where entry $p_{uv}=a_{uv}/k_{u}$ indicates the probability of jumping from $u$ to $v$ in one step. Mathematically, when studying random walks on graph $\mathcal{G}(\mathcal{V},\mathcal{E})$, a significant index for a walker starting out from vertex $u$ is the hitting time $\mathcal{H}_{u\rightarrow v}$ that is in fact the expected time taken by the walker to first reach destination vertex $v$. As a consequence, for the whole graph $\mathcal{G}(\mathcal{V},\mathcal{E})$ in question, mean hitting time $\langle\mathcal{H}_{\mathcal{G}}\rangle$ can be defined as the averaged value over quantities $\mathcal{H}_{u\rightarrow v}$ for all vertex pairs $u$ and $v$, and is given by

\begin{equation}\label{eqa:MF-2-3-1}
\langle\mathcal{H}_{\mathcal{G}}\rangle=\frac{1}{|\mathcal{V}|(|\mathcal{V}|-1)}\sum_{u,v\in\mathcal{V}}\mathcal{H}_{u\rightarrow v}.
\end{equation}

It is well known that for a given graph $\mathcal{G}(\mathcal{V},\mathcal{E})$, one of most commonly used techniques for calculating quantity $\langle\mathcal{H}_{\mathcal{G}}\rangle$ is based on Laplcian matrix \cite{Zhu-1996}. Specifically, this is as follows.

\textbf{Lemma 2.1 \cite{Zhu-1996}} The solution of mean hitting time $\langle\mathcal{H}_{\mathcal{G}}\rangle$ for random walks on  graph $\mathcal{G}(\mathcal{V},\mathcal{E})$ is expressed as

\begin{equation}\label{eqa:MF-2-3-2}
\langle\mathcal{H}_{\mathcal{G}}\rangle=\frac{2|\mathcal{E}|}{|\mathcal{V}|-1}\sum_{i=2}^{|\mathcal{V}|}\frac{1}{\lambda_{i}},
\end{equation}
where $\lambda_{i}$ is all the nonzero eigenvalues of Laplacian matrix $\mathbf{L}_{\mathcal{G}}(=\mathbf{D}_{\mathcal{G}}-\mathbf{A}_{\mathcal{G}})$.

On the other hand, there are some more effective methods for specific graphs, such as, tree, when determining mean hitting time.

\textbf{Lemma 2.2 \cite{Ma-pre-2020}} The solution of mean hitting time $\langle\mathcal{H}_{\mathcal{T}}\rangle$ for random walks on tree $\mathcal{T}$ is given in the following form

\begin{equation}\label{eqa:MF-2-3-3}
\langle\mathcal{H}_{\mathcal{T}}\rangle=\frac{2\mathcal{W}_{\mathcal{T}}}{|\mathcal{T}|},
\end{equation}
in which $|\mathcal{T}|$ represents vertex number of tree $\mathcal{T}$.

As will see later, the closed-form solutions of mean hitting time on trees that are generated below are easily derived using the light shed by Lemma 2.2.

\section{Several primitive operations}

Here, we introduce some primitive operations that have been widely used to create a variety of graph models, such as Vicsek fractal \cite{Vicsek-1983} and T-graph \cite{Redner-2001}. As reported in the rich literature, those models have proven useful in a great number of applications in different fields ranging from network science, graph theory, statistic physics to chemistry, and so forth \cite{Blumen-2004,Gurtovenko-2005}.

\subsection{$m$-order subdivision operation}

For an arbitrary graph $\mathcal{G}(\mathcal{V},\mathcal{E})$, inserting $m$ new vertices into each edge $e_{uv}$ in edge set $\mathcal{E}$ produces a graph $\mathcal{G}^{m}(\mathcal{V}^{m},\mathcal{E}^{m})$, called $m$-order subdivision graph. Such a procedure is often regarded as the \emph{$m$-order subdivision operation}. Equivalently, the end graph $\mathcal{G}^{m}(\mathcal{V}^{m},\mathcal{E}^{m})$ can be also obtained from original graph $\mathcal{G}(\mathcal{V},\mathcal{E})$ by replacing each edge $e_{uv}$ in $\mathcal{E}$ by a path $\mathcal{P}_{uw_{1}\dots w_{m}v}$ of length $(m+1)$ where $w_{i}$ is each newly inserted vertex. Accordingly, a couple of equations associated with $|\mathcal{V}^{m}|$ and $|\mathcal{E}^{m}|$ are given by

 \begin{equation}\label{eqa:MF-3-1-1}
|\mathcal{V}^{m}|=|\mathcal{V}|+m|\mathcal{E}|, \qquad |\mathcal{E}^{m}|=(m+1)|\mathcal{E}|.
\end{equation}
In particular, the case of $m=1$ is in general called edge subdivision in the jargon of graph theory \cite{Bondy-2008}.

\subsection{Type-I growth operation}

For an arbitrary graph $\mathcal{G}(\mathcal{V},\mathcal{E})$, connecting $m$ new vertices as leaves to each vertex $u$ in vertex set $\mathcal{V}$ yields a graph $\mathcal{G}^{\mathrm{I}}_{m}(\mathcal{V}^{\mathrm{I}}_{m},\mathcal{E}^{\mathrm{I}}_{m})$ which, hereafter, we call Type-I graph. Such a procedure is viewed as the \emph{Type-$\mathrm{I}$ growth operation}. As above, we can obtain a system of equations

\begin{equation}\label{eqa:MF-3-2-1}
|\mathcal{V}^{\mathrm{I}}_{m}|=(m+1)|\mathcal{V}|, \qquad |\mathcal{E}^{\mathrm{I}}_{m}|=|\mathcal{E}|+m|\mathcal{V}|.
\end{equation}

\subsection{T-fractal operation}

Given an edge $e_{uv}$, the so-called \emph{$\mathrm{T}$-fractal operation} is described in the following two steps: (1) inserting a new vertex $w$ on edge $e_{uv}$, and (2) connecting an additional vertex $w_{1}$ as leaf to the newly added vertex $w$. More generally, we can have \emph{$m$-order $\mathrm{T}$-fractal operation} if $m$ additional vertices $w_{i}$ ($i\in[1,m]$) are connected to vertex $w$. Consider an arbitrary graph $\mathcal{G}(\mathcal{V},\mathcal{E})$, we can see after applying $m$-order $\mathrm{T}$-fractal operation to each edge $e_{uv}$ in edge set $\mathcal{E}$ that vertex number $|\mathcal{V}^{\mathrm{T}}_{m}|$ and edge number $|\mathcal{V}^{\mathrm{T}}_{m}|$ of the end graph $\mathcal{G}^{\mathrm{T}}_{m}(\mathcal{V}^{\mathrm{T}}_{m},\mathcal{E}^{\mathrm{T}}_{m})$ satisfy

\begin{equation}\label{eqa:MF-3-3-1}
|\mathcal{V}^{\mathrm{T}}_{m}|=|\mathcal{V}|+(m+1)|\mathcal{E}|, \qquad |\mathcal{E}^{\mathrm{T}}_{m}|=(m+2)|\mathcal{E}|.
\end{equation}

\subsection{V-fractal operation}

For a given graph $\mathcal{G}(\mathcal{V},\mathcal{E})$ in which the greatest vertex degree is equal to $k_{max}$, we implement the so-called \emph{$m$-order $\mathrm{V}$-fractal operation} where $m$ is no less than $k_{max}$ as follows: (1) placing two new vertices on each edge $e_{uv}$ in edge set $\mathcal{E}$, and (2) connecting $m-k_{u}$ additional vertices to each vertex $u$ in vertex set $\mathcal{V}$ where $k_{u}$ represents vertex degree. After that, each pre-existing vertex $u$ in graph $\mathcal{G}(\mathcal{V},\mathcal{E})$ is considered $m$-saturated. The resulting graph $\mathcal{G}^{\mathrm{V}}_{m}(\mathcal{V}^{\mathrm{V}}_{m},\mathcal{E}^{\mathrm{V}}_{m})$ follows

\begin{equation}\label{eqa:MF-3-4-1}
|\mathcal{V}^{\mathrm{V}}_{m}|=(m+1)|\mathcal{V}|, \qquad |\mathcal{E}^{\mathrm{V}}_{m}|=|\mathcal{E}|+m|\mathcal{V}|.
\end{equation}

\subsection{Type-II growth operation}

For an arbitrary graph $\mathcal{G}(\mathcal{V},\mathcal{E})$, connecting $mk_{u}$ new vertices as leaves to each vertex $u$ with degree $k_{u}$ in vertex set $\mathcal{V}$ yields a graph $\mathcal{G}^{\mathrm{II}}_{m}(\mathcal{V}^{\mathrm{II}}_{m},\mathcal{E}^{\mathrm{II}}_{m})$. For convenience, the resulting graph $\mathcal{G}^{\mathrm{II}}_{m}(\mathcal{V}^{\mathrm{II}}_{m},\mathcal{E}^{\mathrm{II}}_{m})$ is called $m$-order Type-II graph. Accordingly, such a procedure is thought of as the \emph{$m$-order Type-$\mathrm{II}$ growth operation}. As previously, we can obtain a couple of equations

\begin{equation}\label{eqa:MF-3-5-1}
|\mathcal{V}^{\mathrm{II}}_{m}|=|\mathcal{V}|+2m|\mathcal{E}|, \qquad |\mathcal{E}^{\mathrm{II}}_{m}|=(2m+1)|\mathcal{E}|.
\end{equation}

\subsection{Type-III growth operation}

Consider an arbitrary graph $\mathcal{G}(\mathcal{V},\mathcal{E})$ whose greatest vertex degree is supposed to equal $k_{max}$, we can obtain a growth graph $\mathcal{G}^{\mathrm{III}}_{m}(\mathcal{V}^{\mathrm{III}}_{m},\mathcal{E}^{\mathrm{III}}_{m})$ by attaching $m-k_{u}$ new vertices to each vertex $u$ with degree $k_{u}$. Clearly, parameter $m$ is no less than the greatest degree $k_{max}$. Such a procedure is defined to be \emph{Type-$\mathrm{III}$ growth operation}. As a result, we can write

\begin{equation}\label{eqa:MF-3-6-1}
|\mathcal{V}^{\mathrm{III}}_{m}|=(m+1)|\mathcal{V}|-2|\mathcal{E}|, \qquad |\mathcal{E}^{\mathrm{III}}_{m}|=m|\mathcal{V}|-|\mathcal{E}|.
\end{equation}

\subsection{Other derivatives}

It should be mentioned that the above-defined operations are most fundamental implementations in the current study of networked models \cite{Blumen-2004,Gurtovenko-2005},\cite{Burnashev-2012}-\cite{Bartolo-2009}. And, some more complicated operations can in fact be obtained based on them via various kinds of simple combinatorial manners. As an illustrative example, we do build up an operation by combining $2$-order subdivision operation with $1$-order Type-$\mathrm{II}$ growth operation. Equivalently speaking, we first insert two new vertices on each edge in the initial graph $\mathcal{G}(\mathcal{V},\mathcal{E})$ and then connect $k_{u}$ new vertices to each existing vertex $u$ with degree $k_{u}$ of graph $\mathcal{G}(\mathcal{V},\mathcal{E})$. As shown in the prior work \cite{Ma-arX-2021}, this kind of operations have been used to create a family of networked models with interesting properties including fractal feature. In addition, stochastic versions are also generated via introducing randomness into the process of constructing graphs with respect to primitive operations. Due to the limitation of space, we omit detailed descriptions about development of other derivatives, which is left for interested reader as an exercise.

\emph{Remark 1} As mentioned above, the goal of this work is to determine exact solutions of some structural parameters of various growth trees that are generated using operations proposed above. Therefore, an arbitrary tree $\mathcal{T}$ is always chose to serve as the seed.

\emph{Remark 2} Obviously, given an arbitrary tree $\mathcal{T}$ as seed, each primitive operation mentioned above is implemented iteratively until a desirable model is obtained. For example, we can have a series of growth tree models, which are denoted by $\mathcal{T}^{m}(t),\mathcal{T}^{\mathrm{I}}_{m}(t),\mathcal{T}^{\mathrm{T}}_{m}(t),\mathcal{T}^{\mathrm{V}}_{m}(t),\mathcal{T}^{\mathrm{II}}_{m}(t)$ as well as $\mathcal{T}^{\mathrm{III}}_{m}(t)$ in turn, after $t$ time steps.

\section{Proofs of main results}

This section aims at showing our main results that have been organized in Theorems I.1-I.6. More specifically, the corresponding rigorous proof of each theorem is provided. Note that all the proofs are developed by mean of a more effective manner which we call Mapping Transformation. The thought behind this type of calculational manners is in spirit similar to that from the normalisation group in real space.

In what follows, let us divert our attention to demonstration of main results. First of all, we focus on the simplest case as below.

\subsection{Proof of Theorem I.1}

Given a tree $\mathcal{T}$ as required above, we can abuse $\mathcal{T}$ to denote the corresponding vertex set. Equivalently speaking, we have $|\mathcal{T}|=|\{u:u\in\mathcal{T}\}|=n$. Before beginning our discussions, some necessary notations are listed as below. In view of the concept of $m$-order subdivision operation, it is clear to see that there are $m$ vertices inserted into each edge $e_{uv}$ in tree $\mathcal{T}$. For our purpose, each of these newly inserted vertices into edge $e_{uv}$ is assigned a unique label $w_{e_{uv}}^{i}$ ($i\in[1,m]$), and then they are grouped into a set $\Lambda_{e_{uv}}$, i.e., $\Lambda_{e_{uv}}=\{w_{e_{uv}}^{i}:i\in[1,m]\}$. Based on this, all the new vertices added into tree $\mathcal{T}^{m}(1)$ constitute set $\Lambda^{m}(1)=\bigcup_{e_{uv}\in\mathcal{T}}\Lambda_{e_{uv}}$. To make further progress, we can write $\mathcal{T}^{m}(1)=\mathcal{T}\bigcup\Lambda^{m}(1)$. Now, let us start to validate Theorem I.1 in stages.

\emph{Case 4.1.1} From the concrete construction of tree $\mathcal{T}^{m}(1)$, it is straightforward to see that distance $d'_{uv}$ between vertices $u$ and $v$ in set $\mathcal{T}$ follows $d'_{uv}=(m+1)d_{uv}$ where $d_{uv}$ represents distance of the same pair of vertices in tree $\mathcal{T}$. This further indicates the following expression

\begin{equation}\label{eqa:MF-4-1-1}
\mathcal{W}_{\mathcal{T}^{m}(1)}(1)=\frac{1}{2}\sum_{u,v\in\mathcal{T}}d'_{uv}=(m+1)\mathcal{W}_{\mathcal{T}}.
\end{equation}

\emph{Case 4.1.2} Now, we consider distance $d'_{uw_{e_{xy}}^{i}}$ where vertex $u$ is in set $\mathcal{T}$ and $w_{e_{xy}}^{i}$ belongs to some set $\Lambda_{e_{xy}}$. It is worth noticing that edge $e_{xy}$ may be identical to edge $e_{uv}$. In addition, we require that path $\mathcal{P}_{uy}$ in tree $\mathcal{T}^{m}(1)$ be divided into three segments $\mathcal{P}_{uv},\mathcal{P}_{vx}$ and $\mathcal{P}_{xy}$ when edge $e_{xy}$ is different from edge $e_{uv}$. Next, using a mapping transformation from vertex pair $u$ and $w_{e_{xy}}^{i}$ to vertex pair $u$ and $y$, one finds a close relationship $d'_{uw_{e_{xy}}^{i}}+(m+1-i)=d'_{uy}$. Similarly, there also exists an analog between distances $d'_{w_{e_{uv}}^{j}y}$ and $d'_{uy}$, namely, $d'_{w_{e_{uv}}^{j}y}+(m+1-i)=d'_{uy}$. Taken together, we have

\begin{equation}\label{eqa:MF-4-1-2}
\mathcal{W}_{\mathcal{T}^{m}(1)}(2)=\sum_{u\in\mathcal{T}}\sum_{e_{xy}\in\mathcal{T}}\sum_{w_{e_{xy}}^{i}\in\Lambda_{e_{xy}}}d'_{uw_{e_{xy}}^{i}}=2m\mathcal{W}_{\mathcal{T}^{m}(1)}(1)-2\left(
                                                             \begin{array}{c}
                                                               n \\
                                                               2 \\
                                                             \end{array}
                                                           \right)\left(
                                                             \begin{array}{c}
                                                               m+1 \\
                                                               2 \\
                                                             \end{array}
                                                           \right).
\end{equation}

\emph{Case 4.1.3} The left task is to determine distance between two vertices in set $\Lambda^{m}(1)$. For a pair of vertices $w_{e_{uv}}^{i}$ and $w_{e_{uv}}^{j}$ in an identical sub-set $\Lambda_{e_{uv}}$, we can find after some algebra that the summations over distances $d'_{w_{e_{xy}}^{i}w_{e_{xy}}^{j}}$, denoted by $\mathcal{W}_{\mathcal{T}^{m}(1)}^{1}(3)$, is calculated to equal

\begin{equation}\label{eqa:MF-4-1-3}
\mathcal{W}_{\mathcal{T}^{m}(1)}^{1}(3)=\frac{1}{2}\sum_{w_{e_{uv}}^{i},w_{e_{uv}}^{j}\in\Lambda_{e_{uv}}}d'_{w_{e_{uv}}^{i}w_{e_{uv}}^{j}}=\left(
                                                             \begin{array}{c}
                                                               m+1 \\
                                                               3 \\
                                                             \end{array}
                                                           \right).
\end{equation}
After that, for two vertices $w_{e_{uv}}^{i}$ and $w_{e_{xy}}^{j}$ from distinct sub-sets $\Lambda_{e_{uv}}$ and $\Lambda_{e_{xy}}$, using a similar mapping transformation from vertex pair $w_{e_{uv}}^{i}$ and $w_{e_{xy}}^{j}$ to vertex pair $u$ and $y$ as above leads to a connection of distance $d'_{w_{e_{uv}}^{i}w_{e_{xy}}^{j}}$ to distance $d'_{uy}$, i.e., $d'_{w_{e_{uv}}^{i}w_{e_{xy}}^{j}}+(m+1-i+j)=d'_{uy}$. This certainly suggests the coming expression

\begin{equation}\label{eqa:MF-4-1-4}
\mathcal{W}_{\mathcal{T}^{m}(1)}^{2}(3)=\frac{1}{2}\sum_{w_{e_{uv}}^{i}\in\Lambda_{e_{uv}}}\sum_{w_{e_{xy}}^{j}\in\Lambda_{e_{xy}}}d'_{w_{e_{uv}}^{i}w_{e_{xy}}^{j}}=m^{2}[\mathcal{W}_{\mathcal{T}^{m}(1)}(1)-(n-1)(m+1)]-2m\left(
                                                             \begin{array}{c}
                                                               n-1 \\
                                                               2 \\
                                                             \end{array}
                                                           \right)\left(
                                                             \begin{array}{c}
                                                               m+1 \\
                                                               2 \\
                                                             \end{array}
                                                           \right).
\end{equation}
Using Eqs.(\ref{eqa:MF-4-1-3}) and (\ref{eqa:MF-4-1-4}), the final formula of summation $\mathcal{W}_{\mathcal{T}^{m}(1)}(3)$ over distances of all possible vertex pairs considered herein is written as

\begin{equation}\label{eqa:MF-4-1-5}
\begin{aligned}\mathcal{W}_{\mathcal{T}^{m}(1)}(3)&=(n-1)\mathcal{W}_{\mathcal{T}^{m}(1)}^{1}(3)+\mathcal{W}_{\mathcal{T}^{m}(1)}^{2}(3)\\
&=m^{2}[\mathcal{W}_{\mathcal{T}^{m}(1)}(1)-(n-1)(m+1)]-2m\left(
                                                             \begin{array}{c}
                                                               n-1 \\
                                                               2 \\
                                                             \end{array}
                                                           \right)\left(
                                                             \begin{array}{c}
                                                               m+1 \\
                                                               2 \\
                                                             \end{array}
                                                           \right)+(n-1)\left(
                                                             \begin{array}{c}
                                                               m+1 \\
                                                               3 \\
                                                             \end{array}
                                                           \right)
                                                           \end{aligned}
\end{equation}
And then, substituting Eqs.(\ref{eqa:MF-4-1-1}), (\ref{eqa:MF-4-1-2}) and (\ref{eqa:MF-4-1-5}) into expression $\mathcal{W}_{\mathcal{T}^{m}(1)}=\sum_{i=1}^{3}\mathcal{W}_{\mathcal{T}^{m}(1)}(i)$ yields the same result as in Eq.(\ref{eqa:MF-4-1-0}). To sum up, we complete the proof of Theorem I.1. \qed

\subsection{Proof of Theorem I.2}

As above, we need to take some notations. Based on definition in Section 3.2, there are $m$ new vertices $u_{i}$ ($i\in[1,m]$)
created for each vertex $u$ in tree $\mathcal{T}$. And then, these new vertices $u_{i}$ are contained in set $\Omega_{u}$. As a result, the vertex set of tree $\mathcal{T}^{\mathrm{I}}_{m}(1)$ is composed of two sub-sets $\mathcal{T}$ and $\Omega^{\mathrm{I}}_{m}(1)=\bigcup_{u\in\mathcal{T}}\Omega_{u}$. Similarly, determining Wiener index $\mathcal{W}_{\mathcal{T}^{\mathrm{I}}_{m}(1)}$ is equivalently transformed into calculating three classes of distances as will be shown shortly.

\emph{Case 4.2.1} For a given pair of vertices $u$ and $v$ in sub-set $\mathcal{T}$, Type-I growth operation has no influence on distance $d_{uv}$. This suggests

\begin{equation}\label{eqa:MF-4-2-1}
\mathcal{W}_{\mathcal{T}^{\mathrm{I}}_{m}(1)}(1)=\frac{1}{2}\sum_{u,v\in\mathcal{T}}d'_{uv}=\frac{1}{2}\sum_{u,v\in\mathcal{T}}d_{uv}=\mathcal{W}_{\mathcal{T}}.
\end{equation}

\emph{Case 4.2.2} Next, let us consider distance $d'_{uv_{i}}$ where one vertex is chose from sub-set $\mathcal{T}$ and the other, namely, vertex $v_{i}$, is in some sub-set $\Omega_{v}$. Apparently, the total number of vertex pairs of this kind is equivalent to $n(n+1)m/2$. Among which, there are $n(n-1)m/2$ pairs of vertices in which vertex $v_{i}$ is not in sub-set $\Omega_{u}$. In this case, distance $d'_{uv_{i}}$ is derived by using a mapping transformation as below. It is easy to see that path $\mathcal{P}_{uv_{i}}$ is composed of two segments $\mathcal{P}_{uv}$ and $\mathcal{P}_{vv_{i}}$, which means that distance $d'_{uv_{i}}$ equals $d'_{uv}$ plus $d'_{vv_{i}}$. Therefore, we have

\begin{equation}\label{eqa:MF-4-2-2}
\mathcal{W}_{\mathcal{T}^{\mathrm{I}}_{m}(1)}^{1}(2)=\frac{1}{2}\sum_{u,v(\neq u)\in\mathcal{T}}\sum_{v_{i}\in\Omega_{v}}d'_{uv_{i}}=2m\mathcal{W}_{\mathcal{T}^{\mathrm{I}}_{m}(1)}(1)+2m\left(
                                                             \begin{array}{c}
                                                               n \\
                                                               2 \\
                                                             \end{array}
                                                           \right).
\end{equation}
On the other hand, when two vertices in question are adjacent, distance $d'_{uu_{i}}$ is surely equal to $1$. After that, we obtain the following formula

\begin{equation}\label{eqa:MF-4-2-3}
\mathcal{W}_{\mathcal{T}^{\mathrm{I}}_{m}(1)}^{2}(2)=\sum_{u\in\mathcal{T}}\sum_{u_{i}\in\Omega_{u}}d'_{uu_{i}}=nm.
\end{equation}
Armed with results in Eqs.(\ref{eqa:MF-4-2-2}) and (\ref{eqa:MF-4-2-3}), we write

\begin{equation}\label{eqa:MF-4-2-4}
\mathcal{W}_{\mathcal{T}^{\mathrm{I}}_{m}(1)}(2)=\mathcal{W}_{\mathcal{T}^{\mathrm{I}}_{m}(1)}^{1}(2)+\mathcal{W}_{\mathcal{T}^{\mathrm{I}}_{m}(1)}^{2}(2)=2m\mathcal{W}_{\mathcal{T}^{\mathrm{I}}_{m}(1)}(1)+nm+2m\left(
                                                             \begin{array}{c}
                                                               n \\
                                                               2 \\
                                                             \end{array}
                                                           \right).
\end{equation}

\emph{Case 4.2.3} Lastly, we estimate quantity $\mathcal{W}_{\mathcal{T}^{\mathrm{I}}_{m}(1)}(3)$ that is the summation over distances of all possible pairs of vertices in sub-set $\Omega^{\mathrm{I}}_{m}(1)$. This issue is addressed via considering two cases. The first case is to determine distance $d'_{u_{i}u_{j}}$ in which two vertices are in an identical set $\Omega_{u}$. And then, the following expression is easy to check

\begin{equation}\label{eqa:MF-4-2-5}
\mathcal{W}_{\mathcal{T}^{\mathrm{I}}_{m}(1)}^{1}(3)=\sum_{u\in\mathcal{T}}\sum_{u_{i},u_{j}\in\Omega_{u}}d'_{u_{i}u_{j}}=2n\left(
                                                             \begin{array}{c}
                                                               m \\
                                                               2 \\
                                                             \end{array}
                                                           \right).
\end{equation}
The other case is to evaluate distance $d'_{u_{i}v_{j}}$ where two vertices are from different sets $\Omega_{u}$ and $\Omega_{v}$. It is not hard to find a mapping transformation that path $\mathcal{P}_{u_{i}v_{j}}$ is divided into three segments $\mathcal{P}_{u_{i}u},\mathcal{P}_{uv}$ and $\mathcal{P}_{vv_{j}}$. This further implies

\begin{equation}\label{eqa:MF-4-2-6}
\mathcal{W}_{\mathcal{T}^{\mathrm{I}}_{m}(1)}^{2}(3)=\frac{1}{2}\sum_{u\in\mathcal{T}}\sum_{u_{i}\in\Omega_{u}}\sum_{v(\neq u)\in\mathcal{T}}\sum_{v_{j}\in\Omega_{v}}d'_{u_{i}v_{j}}=m^{2}\mathcal{W}_{\mathcal{T}^{\mathrm{I}}_{m}(1)}(1)+2m^{2}\left(
                                                             \begin{array}{c}
                                                               n \\
                                                               2 \\
                                                             \end{array}
                                                           \right).
\end{equation}
Combining Eq.(\ref{eqa:MF-4-2-5}) with Eq.(\ref{eqa:MF-4-2-6}) yields

\begin{equation}\label{eqa:MF-4-2-7}
\mathcal{W}_{\mathcal{T}^{\mathrm{I}}_{m}(1)}(3)=\mathcal{W}_{\mathcal{T}^{\mathrm{I}}_{m}(1)}^{1}(3)+\mathcal{W}_{\mathcal{T}^{\mathrm{I}}_{m}(1)}^{2}(3)=m^{2}\mathcal{W}_{\mathcal{T}^{\mathrm{I}}_{m}(1)}(1)+2m^{2}\left(
                                                             \begin{array}{c}
                                                               n \\
                                                               2 \\
                                                             \end{array}
                                                           \right)+2n\left(
                                                             \begin{array}{c}
                                                               m \\
                                                               2 \\
                                                             \end{array}
                                                           \right).
\end{equation}

By far, we have enumerated exhaustively all possible cases. Hence, Wiener index $\mathcal{W}_{\mathcal{T}^{\mathrm{I}}_{m}(1)}$ is derived via summing over $\mathcal{W}_{\mathcal{T}^{\mathrm{I}}_{m}(1)}(i)$ ($i\in[1,3]$). After using some simple arithmetics, we complete the proof of Theorem I.2. \qed

\subsection{Proof of Theorem I.3}

In essence, it is obvious to show that $\mathrm{T}$-fractal operation introduced in Section 3.3 is viewed as an extension of subdivision defined in Section 3.1. Therefore, we make use of some previous notations. For example, that vertex inserted into edge $e_{uv}$ is uniquely remarked $w_{e_{uv}}$. And then, all the newly inserted vertices $w_{e_{uv}}$ are collected into set $\Lambda^{\mathrm{T}}_{m}(1)$. In addition, we denote by $\Psi_{w_{e_{uv}}}$ a set consisting of those vertices $w^{i}_{e_{uv}}$ ($i\in[1,m]$) attached to vertex $w_{e_{uv}}$. Accordingly, set $\Psi^{\mathrm{T}}_{m}(1)$ is a collection of sub-sets $\Psi_{w_{e_{uv}}}$, that is to say, $\Psi^{\mathrm{T}}_{m}(1)=\bigcup_{w_{e_{uv}}\in\Lambda^{\mathrm{T}}_{m}(1)}\Psi_{w_{e_{uv}}}$.

Now, the remainder of our tasks are to determine distance associated with each vertex in set $\Psi^{\mathrm{T}}_{m}(1)$. This is due to consequences derived in Theorem I.1. As will see, we proceed other portion of calculations about Wiener index $\mathcal{W}_{\mathcal{T}^{\mathrm{T}}_{m}(1)}$ in stages. Note also that we still adopt a statement that path $\mathcal{P}_{uy}$ in tree $\mathcal{T}^{\mathrm{T}}_{m}(1)$ consists of three segments $\mathcal{P}_{uv},\mathcal{P}_{vx}$ and $\mathcal{P}_{xy}$ when edge $e_{xy}$ is distinct from edge $e_{uv}$.

\emph{Case 4.3.1} For a given pair of vertices $w^{i}_{e_{xy}}$ and $u$, we measure distance $d'_{uw^{i}_{e_{xy}}}$ using a mapping transformation as follows. Path $\mathcal{P}_{uw^{i}_{e_{xy}}}$ in tree $\mathcal{T}^{\mathrm{T}}_{m}(1)$ is obtained from path $\mathcal{P}_{uy}$ by both removing edge $e_{w_{e_{xy}}y}$ and adding an additional edge $e_{w_{e_{xy}}w^{i}_{e_{xy}}}$. This certainly suggests that path $\mathcal{P}_{uw^{i}_{e_{xy}}}$ has a length with path $\mathcal{P}_{uy}$ in common. Also, the latter has been discussed in detail in Case 4.1.1, and, however, parameter $m$ is now equal to $1$. As such, we have

\begin{equation}\label{eqa:MF-4-3-1}
\mathcal{W}_{\mathcal{T}^{\mathrm{I}}_{m}(1)}(1)=\sum_{u\in\mathcal{T}}\sum_{w^{i}_{e_{xy}}\in\Psi^{\mathrm{T}}_{m}(1)}d'_{uw^{i}_{e_{xy}}}=2m\mathcal{W}_{\mathcal{T}^{1}(1)}(1)=4m\mathcal{W}_{\mathcal{T}}.
\end{equation}
Here, we have used Eq.(\ref{eqa:MF-4-1-1}).

\emph{Case 4.3.2} From here on out, let us pay more attention on calculation of distance $d'_{w^{i}_{e_{xy}}w_{e_{uv}}}$ between vertex $w^{i}_{e_{xy}}$ in sub-set $\Psi^{\mathrm{T}}_{m}(1)$ and vertex $w_{e_{uv}}$ in sub-set $\Lambda^{\mathrm{T}}_{m}(1)$. There are in fact two distinct cases. The one is to determine distance $d'_{w^{i}_{e_{uv}}w_{e_{uv}}}$. In this case, quantity $\mathcal{W}_{\mathcal{T}^{\mathrm{I}}_{m}(1)}^{1}(2)$, summation over distances $d'_{w^{i}_{e_{uv}}w_{e_{uv}}}$ of all possible vertex pairs, is calculated to yield

\begin{equation}\label{eqa:MF-4-3-2}
\mathcal{W}_{\mathcal{T}^{\mathrm{I}}_{m}(1)}^{1}(2)=\sum_{w_{e_{uv}}\in\Lambda^{\mathrm{T}}_{m}(1)}\sum_{i=1}^{m}d'_{w^{i}_{e_{uv}}w_{e_{uv}}}=m(n-1).
\end{equation}
The other is to study distances $d'_{w^{i}_{e_{xy}}w_{e_{uv}}}$ where vertex $w^{i}_{e_{xy}}$ is not adjacent to vertex $w_{e_{uv}}$. Here, we employ a mapping transformation as below. Path $\mathcal{P}_{w^{i}_{e_{xy}}w_{e_{uv}}}$ is reduced into path $\mathcal{P}_{w_{e_{xy}}w_{e_{uv}}}$ according to the removal of edge $e_{w^{i}_{e_{xy}}w_{e_{xy}}}$. This further leads to the next equation

\begin{equation}\label{eqa:MF-4-3-3}
\mathcal{W}_{\mathcal{T}^{\mathrm{I}}_{m}(1)}^{2}(2)=\frac{1}{2}\sum_{w_{e_{uv}}\in\Lambda^{\mathrm{T}}_{m}(1)}\sum_{w_{e_{xy}}(\neq w_{e_{uv}})\in\Lambda^{\mathrm{T}}_{m}(1)}\sum_{w^{i}_{e_{xy}}\in\Psi_{w_{e_{xy}}}}d'_{w^{i}_{e_{xy}}w_{e_{uv}}}=2m\mathcal{W}_{\mathcal{T}^{1}(1)}(3)+2m\left(
                                                             \begin{array}{c}
                                                               n-1 \\
                                                               2 \\
                                                             \end{array}
                                                           \right).
\end{equation}
Taken together, we obtain

\begin{equation}\label{eqa:MF-4-3-4}
\mathcal{W}_{\mathcal{T}^{\mathrm{I}}_{m}(1)}(2)=\mathcal{W}_{\mathcal{T}^{\mathrm{I}}_{m}(1)}^{1}(2)+\mathcal{W}_{\mathcal{T}^{\mathrm{I}}_{m}(1)}^{2}(2)=2m\left[2\mathcal{W}_{\mathcal{T}}-2(n-1)-\left(
                                                             \begin{array}{c}
                                                               n-1 \\
                                                               2 \\
                                                             \end{array}
                                                           \right)\right]+m(n-1),
\end{equation}
where we have made use of Eq.(\ref{eqa:MF-4-1-5}).

\emph{Case 4.3.3} Finally, we calculate distance between two vertices selected from set $\Psi^{\mathrm{T}}_{m}(1)$. As above, we encounter two distinct cases, i.e., computing distance $d'_{w^{i}_{e_{xy}}w^{j}_{e_{uv}}}$ and determining quantity $d'_{w^{i}_{e_{uv}}w^{j}_{e_{uv}}}$. First of all, let us focus on the latter case. In an identical sub-set $\Psi_{w_{e_{uv}}}$, distance $d'_{w^{i}_{e_{uv}}w^{j}_{e_{uv}}}$ is equal to $2$ for a pair of vertices $w^{i}_{e_{uv}}$ and $w^{j}_{e_{uv}}$, which means

\begin{equation}\label{eqa:MF-4-3-5}
\mathcal{W}_{\mathcal{T}^{\mathrm{I}}_{m}(1)}^{1}(3)=\frac{1}{2}\sum_{w_{e_{uv}}\in\Lambda^{\mathrm{T}}_{m}(1)}\sum_{w^{i}_{e_{uv}},w^{j}_{e_{uv}}\in\Psi_{w_{e_{uv}}}}d'_{w^{i}_{e_{uv}}w^{j}_{e_{uv}}}=2(n-1)\left(
                                                             \begin{array}{c}
                                                               m \\
                                                               2 \\
                                                             \end{array}
                                                           \right)
\end{equation}
For an arbitrary pari of vertices $w^{i}_{e_{uv}}$ and $w^{j}_{e_{xy}}$ from two distinct sub-sets $\Psi_{w_{e_{uv}}}$ and $\Psi_{w_{e_{xy}}}$, respectively, we will compute distance $d'_{w^{j}_{e_{xy}}w^{i}_{e_{uv}}}$ by utilizing a mapping transformation with respect to distance $d'_{w_{e_{xy}}w_{e_{uv}}}$. More specifically, path $\mathcal{P}_{w^{j}_{e_{xy}}w^{i}_{e_{uv}}}$ can degrade into path $\mathcal{P}_{w_{e_{xy}}w_{e_{uv}}}$ by deleting two end-edges $e_{w^{j}_{e_{xy}}w_{e_{xy}}}$ and $e_{w_{e_{uv}}w^{i}_{e_{uv}}}$. Therefore, we arrive at the coming equation

\begin{equation}\label{eqa:MF-4-3-6}
\begin{aligned}\mathcal{W}_{\mathcal{T}^{\mathrm{I}}_{m}(1)}^{2}(3)&=\frac{1}{2}\sum_{w_{e_{uv}}\in\Lambda^{\mathrm{T}}_{m}(1)}\sum_{w^{i}_{e_{uv}}\in\Psi_{w_{e_{uv}}}}\sum_{w_{e_{xy}}(\neq w_{e_{uv}})\in\Lambda^{\mathrm{T}}_{m}(1)}\sum_{w^{j}_{e_{xy}}\in\Psi_{w_{e_{xy}}}}d'_{w^{j}_{e_{xy}}w^{i}_{e_{uv}}}\\
&=m^{2}\mathcal{W}_{\mathcal{T}^{1}(1)}(3)+2m^{2}\left(
                                                             \begin{array}{c}
                                                               n-1\\
                                                               2 \\
                                                             \end{array}
                                                           \right)
                                                           \end{aligned}.
\end{equation}

From the preceding two equations, quantity $\mathcal{W}_{\mathcal{T}^{\mathrm{I}}_{m}(1)}(3)$ is written as

\begin{equation}\label{eqa:MF-4-3-7}
\mathcal{W}_{\mathcal{T}^{\mathrm{I}}_{m}(1)}(3)=\mathcal{W}_{\mathcal{T}^{\mathrm{I}}_{m}(1)}^{1}(3)+\mathcal{W}_{\mathcal{T}^{\mathrm{I}}_{m}(1)}^{2}(3)=m^{2}\left[2\mathcal{W}_{\mathcal{T}}-2(n-1)\right]+2(n-1)\left(
                                                             \begin{array}{c}
                                                               m \\
                                                               2 \\
                                                             \end{array}
                                                           \right),
\end{equation}
where we perform some fundamental arithmetics and already use Eq.(\ref{eqa:MF-4-1-5}).

In a nutshell, the precise solution of Wiener index $\mathcal{W}_{\mathcal{T}^{\mathrm{T}}_{m}(1)}$ is calculated in the following form

\begin{equation}\label{eqa:MF-4-3-8}
\mathcal{W}_{\mathcal{T}^{\mathrm{T}}_{m}(1)}=\sum_{i=1}^{3}\mathcal{W}_{\mathcal{T}^{\mathrm{I}}_{m}(1)}(i)+\mathcal{W}_{\mathcal{T}^{1}(1)}.
\end{equation}
And, the expression in Eq.(\ref{eqa:MF-4-3-0}) is obtained based on Eq.(\ref{eqa:MF-4-3-8}) by using fundamental calculations and we thus omit it here. This completes the proof of Theorem I.3. \qed

\subsection{Proof of Theorem I.4}

Before starting our discussions, some necessary notations are introduced as below. From the definition of $\mathrm{V}$-fractal operation, there are in fact $m$ new vertices created for each pre-existing vertex $u$ in seed $\mathcal{T}$. Among of them, $k_{u}$ vertices are inserted into all the edges incident with vertex $u$, with each into every edge. The left $(m-k_{u})$ vertices as leaf ones are connected to vertex $u$. Then, these newly added vertices $u_{i}$ adjacent to vertex $u$ are put into the corresponding set $\Gamma_{u}$. After that, all the vertices of tree $\mathcal{T}^{\mathrm{V}}_{m}(1)$ are sorted into two sub-sets $\mathcal{T}$ and $\Gamma_{\mathcal{T}^{\mathrm{V}}_{m}(1)}=\bigcup_{u\in\mathcal{T}}\Gamma_{u}$. So far, we are ready to in detail calculate formula of Winer index $\mathcal{W}_{\mathcal{T}^{\mathrm{V}}_{m}(1)}$. As will be clear to the eye, the following demonstrations are also shown in stages.

\emph{Case 4.4.1} In view of concept defined in Section 3.4, one can see that there are two new vertices placed on each edge $e_{uv}$ in tree $\mathcal{T}$. Such an implementation is in fact $2$-order subdivision. Hence, distance $d'_{uv}$ between vertices $u$ and $v$ in tree $\mathcal{T}^{\mathrm{V}}_{m}(1)$ equals three times larger than distance $d_{uv}$ of the same pair of vertices in seed $\mathcal{T}$. Based on this, we have

\begin{equation}\label{eqa:MF-4-4-1}
\mathcal{W}_{\mathcal{T}^{\mathrm{I}}_{m}(1)}(1)=\frac{1}{2}\sum_{u,v\in\mathcal{T}}d'_{uv}=3\mathcal{W}_{\mathcal{T}}.
\end{equation}

\emph{Case 4.4.2} For a given pair of vertices $u$ and $v_{i}$, we now consider the corresponding distance $d'_{uv_{i}}$. Note that there exist two distinct cases in terms of the selection rule of vertex $v_{i}$. If vertex $v_{i}$ is selected from sub-set $\Gamma_{u}$, in which sense $v_{i}$ is in essence vertex $u_{i}$. And then, we reach the following formula

\begin{equation}\label{eqa:MF-4-4-2}
\mathcal{W}_{\mathcal{T}^{\mathrm{I}}_{m}(1)}^{1}(2)=\sum_{u\in\mathcal{T}}\sum_{u_{i}\in\Gamma_{u}}d'_{uu_{i}}=nm.
\end{equation}
In another case, i.e., vertex $v_{i}$ is not from sub-set $\Gamma_{u}$, we need to use a mapping transformation to quantify distance $d'_{uv_{i}}$ as will be stated shortly. It is known that there are $m$ paths $\mathcal{P}_{uv_{i}}$ in both vertex $u$ and sub-set $\Gamma_{v}$. Among which, $(m-1)$ paths $\mathcal{P}_{uv_{i}}$ is obtained from a unique path $\mathcal{P}_{uv}$ by adding an end-edge $e_{v_{i}v}$. However, the left one path $\mathcal{P}_{uv'_{i}}$ is generated based on path $\mathcal{P}_{uv}$ via deleting an additional edge $e_{v'_{i}v}$. At the same time, the above-mentioned statement holds for both vertex $v$ and sub-set $\Gamma_{u}$ by symmetry. After that, we can conclude

\begin{equation}\label{eqa:MF-4-4-3}
\mathcal{W}_{\mathcal{T}^{\mathrm{I}}_{m}(1)}^{2}(2)=\frac{1}{2}\sum_{u\in\mathcal{T}}\sum_{v(\neq u)\in\mathcal{T}}\sum_{v_{i}\in\Gamma_{v}}d'_{uv_{i}}=2m\mathcal{W}_{\mathcal{T}^{\mathrm{I}}_{m}(1)}(1)+2(m-2)\left(
                                                             \begin{array}{c}
                                                               n \\
                                                               2 \\
                                                             \end{array}
                                                           \right).
\end{equation}
To make further progress, quantity $\mathcal{W}_{\mathcal{T}^{\mathrm{I}}_{m}(1)}(2)$ is

\begin{equation}\label{eqa:MF-4-4-4}
\mathcal{W}_{\mathcal{T}^{\mathrm{I}}_{m}(1)}(2)=\mathcal{W}_{\mathcal{T}^{\mathrm{I}}_{m}(1)}^{1}(2)+\mathcal{W}_{\mathcal{T}^{\mathrm{I}}_{m}(1)}^{2}(2)=2m\mathcal{W}_{\mathcal{T}^{\mathrm{I}}_{m}(1)}(1)+nm+2(m-2)\left(
                                                             \begin{array}{c}
                                                               n \\
                                                               2 \\
                                                             \end{array}
                                                           \right).
\end{equation}

\emph{Case 4.4.3} Along the same demonstration as above, let us determine distance $d'_{u_{i}v_{j}}$ where vertices $u_{i}$ and $v_{j}$ are now from set $\Gamma_{\mathcal{T}^{\mathrm{V}}_{m}(1)}$. Analogously, it should be pointed out that two types of cases need be analyzed in detail. When this pair of vertices belong to an identical sub-set $\Gamma_{u}$, its corresponding distance $d'_{u_{i}u_{j}}$ is clearly equivalent to $2$. And then, we can obtain

\begin{equation}\label{eqa:MF-4-4-5}
\mathcal{W}_{\mathcal{T}^{\mathrm{I}}_{m}(1)}^{1}(3)=\frac{1}{2}\sum_{u\in\mathcal{T}}\sum_{u_{i},u_{j}\in\Gamma_{u}}d'_{u_{i}u_{j}}=2n\left(
                                                             \begin{array}{c}
                                                               m \\
                                                               2 \\
                                                             \end{array}
                                                           \right).
\end{equation}
On the other hand, if two vertices $u_{i}$ and $v_{j}$ come from sub-sets $\Gamma_{u}$ and $\Gamma_{v}$, respectively, then the corresponding distance $d'_{u_{i}v_{j}}$ is measured by building up a mapping transformation to quantity $d'_{uv}$. To put this another way, we need to use quantity $d'_{uv}$ as an intermediate variable in the process of determining distance $d'_{u_{i}v_{j}}$. From the description about $\mathrm{V}$-fractal operation, one can understand that there must be $(m-1)^{2}$ paths $\mathcal{P}_{u_{_{i}}v_{j}}$ in both sub-sets $\Gamma_{u}$ and $\Gamma_{v}$ generated from path $\mathcal{P}_{uv}$ by adding two edges $e_{uu_{i}}$ and $e_{vv_{j}}$. In the meantime, there must be a unique path $\mathcal{P}_{u'_{_{i}}v'_{j}}$ obtained based on path $\mathcal{P}_{uv}$ by removing two end-edges $e_{uu'_{i}}$ and $e_{vv'_{j}}$. Additionally, the left $2(m-1)$ paths are all created by means of both the removal of an end-edge in path $\mathcal{P}_{uv}$ and the addition of an additional edge into path $\mathcal{P}_{uv}$. For instance, path $\mathcal{P}_{u_{_{i}}v'_{j}}$ is reduced into $\mathcal{P}_{uv}$ through deleting end-edge $e_{u_{_{i}}u}$ and adding new edge $e_{v'_{j}v}$ simultaneously. Keep this in mind, we obtain

\begin{equation}\label{eqa:MF-4-4-6}
\mathcal{W}_{\mathcal{T}^{\mathrm{I}}_{m}(1)}^{2}(3)=\frac{1}{2}\sum_{u\in\mathcal{T}}\sum_{v(\neq u)\in\mathcal{T}}\sum_{u_{i}\in\Gamma_{u}}\sum_{v_{j}\in\Gamma_{u}}d'_{u_{i}v_{j}}=m^{2}\mathcal{W}_{\mathcal{T}^{\mathrm{I}}_{m}(1)}(1)+2[(m-1)^{2}-1]\left(
                                                             \begin{array}{c}
                                                               n \\
                                                               2 \\
                                                             \end{array}
                                                           \right).
\end{equation}
As a consequence, we have

\begin{equation}\label{eqa:MF-4-4-7}
\mathcal{W}_{\mathcal{T}^{\mathrm{I}}_{m}(1)}(3)=\mathcal{W}_{\mathcal{T}^{\mathrm{I}}_{m}(1)}^{1}(3)+\mathcal{W}_{\mathcal{T}^{\mathrm{I}}_{m}(1)}^{2}(3)=m^{2}\mathcal{W}_{\mathcal{T}^{\mathrm{I}}_{m}(1)}(1)+2n\left(
                                                             \begin{array}{c}
                                                               m \\
                                                               2 \\
                                                             \end{array}
                                                           \right)+2[(m-1)^{2}-1]\left(
                                                             \begin{array}{c}
                                                               n \\
                                                               2 \\
                                                             \end{array}
                                                           \right).
\end{equation}

In order to validate Eq.(\ref{eqa:MF-4-4-0}), the last step is to sum over $\mathcal{W}_{\mathcal{T}^{\mathrm{I}}_{m}(1)}(i)$ ($i\in[1,3]$) and then to perform some fundamental arithmetics. So, we omit the detailed calculations here. In a word, this completes the proof of Theorem I.4. \qed

\subsection{Proof of Theorem I.5}

Similarly, we first need to take some notations used later. According to definition of Type-II growth operation, there will be $mk_{u}$ vertices attached to each degree $k_{u}$ vertex $u$ in seed $\mathcal{T}$. For brevity, these newly created vertices $u_{i}$ ($i\in[1,mk_{u}]$) for pre-existing vertex $u$ are clustered into set $\Theta_{u}$. As such, all the vertices of tree $\mathcal{T}^{\mathrm{II}}_{m}(1)$ are classified into two sub-sets $\mathcal{T}$ and $\Theta_{\mathcal{T}^{\mathrm{II}}_{m}(1)}=\bigcup_{u\in\mathcal{T}}\Theta_{u}$. Bear it in mind, we will deal with calculation of Wiener index $\mathcal{W}_{\mathcal{T}^{\mathrm{II}}_{m}(1)}$ in stages as follows.

\emph{Case 4.5.1} Obviously, Type-II growth operation has no influence on calculation of distance $d'_{uv}$ between two vertices $u$ and $v$ in sub-set $\mathcal{T}$. In other words, distance $d'_{uv}$ is equal to quantity $d_{uv}$, which leads to the following equality

\begin{equation}\label{eqa:MF-4-5-1}
\mathcal{W}_{\mathcal{T}^{\mathrm{II}}_{m}(1)}(1)=\frac{1}{2}\sum_{u,v\in\mathcal{T}}d'_{uv}=\mathcal{W}_{\mathcal{T}}.
\end{equation}

\emph{Case 4.5.2} Below, we focus mainly on calculation of distance between a pair of vertices from sub-sets $\mathcal{T}$ and $\Theta_{\mathcal{T}^{\mathrm{II}}_{m}(1)}$, respectively. Beginning with the simplest case of $d'_{uu_{i}}$, we can straightforwardly write an equality in the following form

\begin{equation}\label{eqa:MF-4-5-2}
\mathcal{W}_{\mathcal{T}^{\mathrm{II}}_{m}(1)}^{1}(2)=\sum_{u\in\mathcal{T}}\sum_{u_{i}\in \Theta_{u}}d'_{uu_{i}}=m\sum_{u\in\mathcal{T}}k_{u}=2m(n-1).
\end{equation}
Following the above analysis, we consider the other case in which vertex $v_{i}$ is from some sub-set $\Theta_{v}$ distinct from $\Theta_{u}$. Using a mapping transformation from path $\mathcal{P}_{uv}$ to path $\mathcal{P}_{uv_{i}}$, one can find that distance $d'_{uv_{i}}$ is equal to distance $d'_{uv}$ plus one. This further indicates

\begin{equation}\label{eqa:MF-4-5-3}
\mathcal{W}_{\mathcal{T}^{\mathrm{II}}_{m}(1)}^{2}(2)=\frac{1}{2}\sum_{u\in\mathcal{T}}\sum_{v(\neq u)\in\mathcal{T}}\sum_{v_{i}\in \Theta_{v}}d'_{uv_{i}}=\frac{m}{2}\sum_{u,v\in\mathcal{T}}(k_{u}+k_{v})d'_{uv}+m(n-1)\sum_{u\in\mathcal{T}}k_{u}.
\end{equation}
It should be mentioned that the first term on the right-hand side of the second equality is not easy to calculate in its present form. To address this issue, we need to appeal to a connection of $2$-order subdivision operation to a specific case of Type-II growth operation where parameter $m$ is supposed to equal $1$, which is in depth explained as below. In which case, only $k_{u}$ vertices $u_{i}$ ($i\in[1,k_{u}]$) are connected to pre-existing vertex $u$ with degree $k_{u}$ in seed $\mathcal{T}$. As a consequence, there are two new vertices created for each existing edge $e_{uv}$ in tree $\mathcal{T}$. Note also that the both new vertices serve as leaf ones in form. On the other hand, one can find from the concept of $2$-order subdivision operation that there are also two new vertices created for each existing edge $e_{uv}$ in tree $\mathcal{T}$. As opposed to Type-II growth operation, they are now inserted into edge $e_{uv}$. Without loss of generality, we assign two labels $u^{\star}_{i}$ and $v^{\star}_{j}$ to the both vertices newly inserted on edge $e_{uv}$. As such, it is clear to see from tree $\mathcal{T}^{2}(1)$ that each degree $k_{u}$ vertex $u$ in set $\mathcal{T}$ is connected to $k_{u}$ new vertices $u^{\star}_{i}$ ($i\in[1,k_{u}]$). And then, there is certainly a one-to-one mapping between vertex $u_{i}$ and $u^{\star}_{i}$. As will be shown above, an analog of quantity $\frac{1}{2}\sum_{u,v\in\mathcal{T}}(k_{u}+k_{v})d'_{uv}$ has been analytically reported in Eq.(\ref{eqa:MF-4-1-2}). So, using both result in Eq.(\ref{eqa:MF-4-1-2}) and the above analysis yields

\begin{equation}\label{eqa:MF-4-5-4}
\frac{1}{2}\sum_{u,v\in\mathcal{T}}(k_{u}+k_{v})d'_{uv}=\frac{1}{3}\left[\mathcal{W}_{\mathcal{T}^{2}(1)}(2)-\mathcal{W}_{\mathcal{T}^{\mathrm{II}}_{1}(1)}^{1}(2)-(n-1)\sum_{u\in\mathcal{T}}k_{u}+2\times2\left(
                                                             \begin{array}{c}
                                                               n \\
                                                               2 \\
                                                             \end{array}
                                                           \right)\right].
\end{equation}
Here, the factor of $1/3$ on the right-hand side of the preceding equation is explained as follows: for a given pair of vertices $u$ and $v$ in set $\mathcal{T}$, the corresponding distance $d'_{uv}$ in tree $\mathcal{T}^{2}(1)$ is three times larger than distance $d_{uv}$ between the same pair of vertices in tree $\mathcal{T}$. The first factor of $2$ in the square bracket on the right-hand side of Eq.(\ref{eqa:MF-4-5-4}) is due to equality $(d'_{u^{\star}_{i}v}+1)/3=d_{uv}=d'_{u_{i}v}-1$ $(\text{or}\; d'_{u^{\star}_{i}v}=3d'_{u_{i}v}-4)$ when considering both a given pair of vertices $u_{i}$ and $v$ in tree $\mathcal{T}^{\mathrm{II}}_{1}(1)$ and the accompanying vertex pair  $u^{\star}_{i}$ and $v$ in tree $\mathcal{T}^{2}(1)$. In a nutshell, we can obtain

\begin{equation}\label{eqa:MF-4-5-5}
\mathcal{W}_{\mathcal{T}^{\mathrm{II}}_{m}(1)}(2)=\mathcal{W}_{\mathcal{T}^{\mathrm{II}}_{m}(1)}^{1}(2)+\mathcal{W}_{\mathcal{T}^{\mathrm{II}}_{m}(1)}^{2}(2)=4m\mathcal{W}_{\mathcal{T}}+\frac{4m}{3}n(n-1)-\frac{2m}{3}\left(
                                                             \begin{array}{c}
                                                               n \\
                                                               2 \\
                                                             \end{array}
                                                           \right)
\end{equation}

\emph{Case 4.5.3} Now, we study distance between an arbitrary pair of vertices chose from sub-set $\Theta_{\mathcal{T}^{\mathrm{II}}_{m}(1)}$. As previously, there are also two distinct cases in such a situation. The first is to estimate distance $d'_{u_{i}u_{j}}$ of two vertices $u_{i}$ and $u_{j}$ in an identical sub-set $\Theta_{u}$. It is easy to estimate, and we can write

\begin{equation}\label{eqa:MF-4-5-6}
\mathcal{W}_{\mathcal{T}^{\mathrm{II}}_{m}(1)}^{1}(3)=\frac{1}{2}\sum_{u\in\mathcal{T}}\sum_{u_{i},u_{j}\in\Theta_{u}}d'_{u_{i}u_{j}}=2\sum_{u\in\mathcal{T}}\left(
                                                             \begin{array}{c}
                                                               mk_{u} \\
                                                               2 \\
                                                             \end{array}
                                                           \right).
\end{equation}
It is worth noticing that determining the last term of the preceding equation directly is a tough task. Thus, we need to perform a simple transformation as follows.

\begin{equation}\label{eqa:MF-4-5-7}
2\sum_{u\in\mathcal{T}}\left(
                                                             \begin{array}{c}
                                                               mk_{u} \\
                                                               2 \\
                                                             \end{array}
                                                           \right)=m^{2}\sum_{u}k_{u}^{2}-m\sum_{u\in\mathcal{T}}k_{u}.
\end{equation}
As will see, the first term on the right-hand side in Eq.(\ref{eqa:MF-4-5-7}) can be derived from the coming calculations reported in Eq.(\ref{eqa:MF-4-5-9}).

Next, when two vertices $u_{i}$ and $v_{j}$ are from distinct sub-sets $\Theta_{u}$ and $\Theta_{v}$, respectively, the corresponding distance $d'_{u_{i}v_{j}}$ can be by definition written as

\begin{equation}\label{eqa:MF-4-5-8}
\mathcal{W}_{\mathcal{T}^{\mathrm{II}}_{m}(1)}^{2}(3)=\frac{1}{2}\sum_{u\in\mathcal{T}}\sum_{u_{i}\in\Theta_{u}}\sum_{v(\neq u)\in\mathcal{T}}\sum_{v_{j}\in\Theta_{v}}d'_{u_{i}v_{j}}=\frac{m^{2}}{2}\sum_{u,v\in\mathcal{T}}k_{u}k_{v}d'_{uv}+m^{2}\sum_{u,v\in\mathcal{T}}k_{u}k_{v}.
\end{equation}
Here, we only need to determine the first term on the right-hand side in the last equality. This is because the second term is easily determined using the following formula

\begin{equation}\label{eqa:MF-4-5-9}
m^{2}\sum_{u,v\in\mathcal{T}}k_{u}k_{v}+m^{2}\sum_{u}k_{u}^{2}=\left(\sum_{u}mk_{u}\right)^{2}.
\end{equation}
For determination of the first term, we still appeal to that connection stated in Case 4.5.2. Specifically, we must make use of the result in Eq.(\ref{eqa:MF-4-1-4}). Due to a similar explanation as in Case 4.5.2, we omit detailed calculations and straightforwardly give the following expression for readability.

\begin{equation}\label{eqa:MF-4-5-10}
\frac{1}{2}\sum_{u,v\in\mathcal{T}}k_{u}k_{v}d'_{uv}=\frac{1}{3}\left\{\mathcal{W}_{\mathcal{T}^{2}(1)}(3)-\mathcal{W}_{\mathcal{T}^{\mathrm{II}}_{1}(1)}^{1}(3)-\sum_{u,v\in\mathcal{T}}k_{u}k_{v}+4\left[\left(
                                                             \begin{array}{c}
                                                               n \\
                                                               2 \\
                                                             \end{array}
                                                           \right)+\left(
                                                             \begin{array}{c}
                                                               n-1 \\
                                                               2 \\
                                                             \end{array}
                                                           \right)\right]\right\}.
\end{equation}
Combining Eq.(\ref{eqa:MF-4-5-6}) with Eq.(\ref{eqa:MF-4-5-8}) produces

\begin{equation}\label{eqa:MF-4-5-11}
\begin{aligned}\mathcal{W}&_{\mathcal{T}^{\mathrm{II}}_{m}(1)}(3)=\mathcal{W}_{\mathcal{T}^{\mathrm{II}}_{m}(1)}^{1}(3)+\mathcal{W}_{\mathcal{T}^{\mathrm{II}}_{m}(1)}^{2}(3)\\
&=4m^{2}\mathcal{W}_{\mathcal{T}}-(n-1)(2n-1)m^{2}+\frac{8m^{2}}{3}(n-1)^{2}-2m(n-1)
+\frac{4m^{2}}{3}\left[\left(
                                                             \begin{array}{c}
                                                               n \\
                                                               2 \\
                                                             \end{array}
                                                           \right)+\left(
                                                             \begin{array}{c}
                                                               n-1 \\
                                                               2 \\
                                                             \end{array}
                                                           \right)\right]
\end{aligned}.
\end{equation}
Here, we have used Eqs.(\ref{eqa:MF-4-1-5}), (\ref{eqa:MF-4-5-6}), (\ref{eqa:MF-4-5-7}) and (\ref{eqa:MF-4-5-9}).

Taken together, substituting Eqs.(\ref{eqa:MF-4-5-1}),(\ref{eqa:MF-4-5-5}) and (\ref{eqa:MF-4-5-11}) into the following summation

\begin{equation}\label{eqa:MF-4-5-12}
\mathcal{W}_{\mathcal{T}^{\mathrm{II}}_{m}(1)}=\sum_{i=1}^{3}\mathcal{W}_{\mathcal{T}^{\mathrm{II}}_{m}(1)}(i),
\end{equation}
and performing some fundamental arithmetics together yields the same result as shown in Eq.(\ref{eqa:MF-4-5-0}), which implies that we finish the proof of Theorem I.5. It should be pointed out that another proof of Theorem I.5 has been reported in the prior work \cite{Ma-arX-2021}. Interested reader is encouraged to refer to \cite{Ma-arX-2021} for more details. \qed

\subsection{Proof of Theorem I.6}

Now, we come to the last theorem in this work. To validate Theorem I.6, we also need to introduce some notations as follows. By definition of Type-III growth operation, one can see that there are $(m-k_{u})$ new vertices $u_{i}$ connected to existing vertex $u$ with degree $k_{u}$ in seed $\mathcal{T}$. For convenience, we use symbol $\Upsilon_{u}$ to represent set consisting of new vertices $u_{i}$. As a result, all the vertices of tree $\mathcal{T}^{\mathrm{III}}_{m}(1)$ are sorted into two mutually disjoint sub-sets $\mathcal{T}$ and $\Upsilon_{\mathcal{T}^{\mathrm{I}}_{m}(1)}=\bigcup_{u\in\mathcal{T}}\Upsilon_{u}$. Then, let us start to demonstrate a detailed analysis in stages.

\emph{Case 4.6.1} For an arbitrary pair of vertices $u$ and $v$ in sub-set $\mathcal{T}$, it is straightforward to obtain the next formula

\begin{equation}\label{eqa:MF-4-6-1}
\mathcal{W}_{\mathcal{T}^{\mathrm{III}}_{m}(1)}(1)=\frac{1}{2}\sum_{u,v\in\mathcal{T}}d'_{uv}=\mathcal{W}_{\mathcal{T}}.
\end{equation}

\emph{Case 4.6.2} Let us turn our attention on calculation of distance between two vertices where one vertex is in sub-set $\mathcal{T}$ and the other is selected from sub-set $\Upsilon_{\mathcal{T}^{\mathrm{I}}_{m}(1)}$. Among of them, we first study this type of distance $d'_{uu_{i}}$ whose two vertices are adjacent to one another. In this situation, we can without difficulty write

\begin{equation}\label{eqa:MF-4-6-2}
\mathcal{W}_{\mathcal{T}^{\mathrm{III}}_{m}(1)}^{1}(2)=\sum_{u\in\mathcal{T}}\sum_{u_{i}\in\Upsilon_{u}}d'_{uu_{i}}=\sum_{u\in\mathcal{T}}(m-k_{u})=mn-2(n-1),
\end{equation}
in which we have taken advantage of Eqs.(\ref{eqa:MF-4-2-3}) and (\ref{eqa:MF-4-5-2}). Next, we consider distance $d'_{uv_{i}}$ in which vertex $v_{i}$ is now in some sub-set $\Upsilon_{v}$ different from $\Upsilon_{u}$ . To put it another way, vertex $u$ is no longer connected to vertex $v_{i}$ by an edge. Facing with this case, one can by definition have

\begin{equation}\label{eqa:MF-4-6-3}
\mathcal{W}_{\mathcal{T}^{\mathrm{III}}_{m}(1)}^{2}(2)=\frac{1}{2}\sum_{u\in\mathcal{T}}\sum_{v(\neq u)\in\mathcal{T}}\sum_{v_{i}\in \Upsilon_{v}}d'_{uv_{i}}=\frac{1}{2}\sum_{u,v\in\mathcal{T}}[(m-k_{u})+(m-k_{v})]d'_{uv}+\sum_{u,v\in\mathcal{T}}[(m-k_{u})+(m-k_{v})].
\end{equation}
Using Eq.(\ref{eqa:MF-4-5-3}), the previous equation can be rearranged as

\begin{equation}\label{eqa:MF-4-6-4}
\begin{aligned}\mathcal{W}_{\mathcal{T}^{\mathrm{III}}_{m}(1)}^{2}(2)&=2m\mathcal{W}_{\mathcal{T}}+mn(n-1)-\mathcal{W}_{\mathcal{T}^{\mathrm{II}}_{1}(1)}^{2}(2)\\
&=2(m-2)\mathcal{W}_{\mathcal{T}}+n(n-1)\left(m-\frac{4}{3}\right)+2(n-1)+\frac{2}{3}\left(
                                                             \begin{array}{c}
                                                               n \\
                                                               2 \\
                                                             \end{array}
                                                           \right)
\end{aligned}.
\end{equation}
We have made use of Eqs.(\ref{eqa:MF-4-5-2}) and (\ref{eqa:MF-4-5-5}). After that, quantity $\mathcal{W}_{\mathcal{T}^{\mathrm{III}}_{m}(1)}(2)$ is shown in the following form

\begin{equation}\label{eqa:MF-4-6-5}
\mathcal{W}_{\mathcal{T}^{\mathrm{III}}_{m}(1)}(2)=\mathcal{W}_{\mathcal{T}^{\mathrm{III}}_{m}(1)}^{1}(2)+\mathcal{W}_{\mathcal{T}^{\mathrm{III}}_{m}(1)}^{2}(2)=2(m-2)\mathcal{W}_{\mathcal{T}}+n(n-1)\left(m-\frac{4}{3}\right)+mn+\frac{2}{3}\left(
                                                             \begin{array}{c}
                                                               n \\
                                                               2 \\
                                                             \end{array}
                                                           \right).
\end{equation}

\emph{Case 4.6.3} In the sequel, we focus on distance of a couple of vertices from sub-set $\Upsilon_{\mathcal{T}^{\mathrm{I}}_{m}(1)}$. As previously, two distinct cases need to be analyzed. First of all, when this pair of vertices are chose from an identical sub-set $\Upsilon_{u}$, it is clear to see the next expression

\begin{equation}\label{eqa:MF-4-6-6}
\mathcal{W}_{\mathcal{T}^{\mathrm{III}}_{m}(1)}^{1}(3)=\frac{1}{2}\sum_{u\in\mathcal{T}}\sum_{u_{i},u_{j}\in\Upsilon_{u}}d'_{u_{i}u_{j}}=2\sum_{u\in\mathcal{T}}\left(
                                                             \begin{array}{c}
                                                               m-k_{u} \\
                                                               2 \\
                                                             \end{array}
                                                           \right)=m(m-1)n-2(2m-1)(n-1)+\sum_{u\in\mathcal{T}}k^{2}_{u}.
\end{equation}
On the other hand, for distance $d'_{u_{i}v_{j}}$ where vertices $u_{i}$ and $v_{j}$ are selected from two distinct sub-sets $\Upsilon_{u}$ and $\Upsilon_{v}$, respectively, one can write

\begin{equation}\label{eqa:MF-4-6-7}
\begin{aligned}\mathcal{W}&_{\mathcal{T}^{\mathrm{III}}_{m}(1)}^{2}(3)=\frac{1}{2}\sum_{u\in\mathcal{T}}\sum_{u_{i}\in\Upsilon_{u}}\sum_{v(\neq u)\in\mathcal{T}}\sum_{v_{j}\in\Upsilon_{v}}d'_{u_{i}v_{j}}=\frac{1}{2}\sum_{u,v\in\mathcal{T}}(m-k_{u})(m-k_{v})d'_{uv}+\sum_{u,v\in\mathcal{T}}(m-k_{u})(m-k_{v})\\
&=m^{2}\mathcal{W}_{\mathcal{T}}-\frac{m}{2}\sum_{u,v\in\mathcal{T}}(k_{u}+k_{v})d'_{uv}-m\sum_{u,v\in\mathcal{T}}(k_{u}+k_{v})+m^{2}n(n-1)+\frac{1}{2}\sum_{u,v\in\mathcal{T}}k_{u}k_{v}d'_{uv}+\sum_{u,v\in\mathcal{T}}k_{u}k_{v}
\end{aligned}
\end{equation}
With Eqs.(\ref{eqa:MF-4-5-3}) and (\ref{eqa:MF-4-5-8}), Eq.(\ref{eqa:MF-4-6-7}) is reorganized as

\begin{equation}\label{eqa:MF-4-6-8}
\mathcal{W}_{\mathcal{T}^{\mathrm{III}}_{m}(1)}^{2}(3)=m^{2}\mathcal{W}_{\mathcal{T}}+m^{2}n(n-1)-2m(n-1)^{2}+\mathcal{W}_{\mathcal{T}^{\mathrm{II}}_{1}(1)}^{2}(3)-\mathcal{W}_{\mathcal{T}^{\mathrm{II}}_{m}(1)}^{2}(2).
\end{equation}
To make further progress, we arrive at quantity $\mathcal{W}_{\mathcal{T}^{\mathrm{III}}_{m}(1)}(3)$ as below

\begin{equation}\label{eqa:MF-4-6-9}
\begin{aligned}\mathcal{W}_{\mathcal{T}^{\mathrm{III}}_{m}(1)}(3)&=\mathcal{W}_{\mathcal{T}^{\mathrm{III}}_{m}(1)}^{1}(3)+\mathcal{W}_{\mathcal{T}^{\mathrm{III}}_{m}(1)}^{2}(3)\\
&=(m-2)^{2}\mathcal{W}_{\mathcal{T}}+n(n-1)\left(m^{2}-\frac{4m}{3}\right)+m(m-1)n-2(n-1)(m-1)\\
&\quad-(n-1)(2n-1)+\left(\frac{8}{3}-2m\right)(n-1)^{2}+\frac{2m}{3}\left(
                                                             \begin{array}{c}
                                                               n \\
                                                               2 \\
                                                             \end{array}
                                                           \right)+\frac{4}{3}\left[\left(
                                                             \begin{array}{c}
                                                               n \\
                                                               2 \\
                                                             \end{array}
                                                           \right)+\left(
                                                             \begin{array}{c}
                                                               n-1 \\
                                                               2 \\
                                                             \end{array}
                                                           \right)\right]
\end{aligned},
\end{equation}
in which Eqs.(\ref{eqa:MF-4-5-2}),(\ref{eqa:MF-4-5-5}),(\ref{eqa:MF-4-5-7}) and (\ref{eqa:MF-4-5-11}) have been employed.

Until now, we have enumerated all possible cases exhaustively. Accordingly, the result in Eq.(\ref{eqa:MF-4-6-0}) is derived based on summation $\mathcal{W}_{\mathcal{T}^{\mathrm{III}}_{m}(1)}=\sum_{i=1}^{3}\mathcal{W}_{\mathcal{T}^{\mathrm{III}}_{m}(1)}(i)$ after some elementary arithmetics. Thus, we finish the proof of Theorem I.6. \qed

In a word, we have proven Theorems I.1-I.6 in a mathematically rigorous manner. It is worth mentioning that all the proofs are developed based on Mapping Transformation established by us. At the same time, we would like to stress that although Wiener indices on six kinds of trees are derived analytically, the corresponding representations seem slightly complicated. In principle, they can look elegant after performing some fundamental arithmetics. Towards this end, we provide the simplified versions below. Accordingly, a concise form also enables us to derive many other interesting structural parameters as will be stated shortly.

\emph{Remark 3} The simplified versions of exact formulas for Wiener indices of trees $\mathcal{T}^{m}(t)$, $\mathcal{T}^{\mathrm{I}}_{m}(t)$, $\mathcal{T}^{\mathrm{T}}_{m}(t)$, $\mathcal{T}^{\mathrm{V}}_{m}(t)$, $\mathcal{T}^{\mathrm{II}}_{m}(t)$ as well as $\mathcal{T}^{\mathrm{III}}_{m}(t)$ are given by, respectively,

\begin{subequations}
\label{eq:whole}
\begin{eqnarray}
\mathcal{W}_{\mathcal{T}^{m}(1)}=(m+1)^{3}\mathcal{W}_{\mathcal{T}}-\frac{m(m+1)^{2}}{2}n^{2}+\frac{m(m+1)(2m+1)}{3}n-\frac{m(m^{2}-1)}{6},\label{subeq:MF-4-7-0-1}
\end{eqnarray}
\begin{equation}
\mathcal{W}_{\mathcal{T}^{\mathrm{I}}_{m}(1)}=(m+1)^{2}\mathcal{W}_{\mathcal{T}}+m(m+1)n^{2}-mn,\label{subeq:MF-4-7-0-2}
\end{equation}
\begin{equation}
\mathcal{W}_{\mathcal{T}^{\mathrm{T}}_{m}(1)}=2(m+2)^{2}\mathcal{W}_{\mathcal{T}}-(m+2)n^{2}-(m-1)(m+2)n+m^{2}+2m,\label{subeq:MF-4-7-0-3}
\end{equation}
\begin{equation}
\mathcal{W}_{\mathcal{T}^{\mathrm{V}}_{m}(1)}=3(m+1)^{2}\mathcal{W}_{\mathcal{T}}+(m-2)(m+1)n^{2}+(m+2)n,\label{subeq:MF-4-7-0-4}
\end{equation}
\begin{equation}
\mathcal{W}_{\mathcal{T}^{\mathrm{II}}_{m}(1)}=(2m+1)^{2}\mathcal{W}_{\mathcal{T}}+m(2m+1)n^{2}-m(5m+3)n+m(3m+2),\label{subeq:MF-4-7-0-5}
\end{equation}
\begin{equation}
\mathcal{W}_{\mathcal{T}^{\mathrm{III}}_{m}(1)}=(m-1)^{2}\mathcal{W}_{\mathcal{T}}+(m-1)^{2}n^{2}+2(m-1)n+1.\label{subeq:MF-4-7-0-6}
\end{equation}
\end{subequations}

\emph{Remark 4} From Eqs.(\ref{subeq:MF-4-7-0-1})-(\ref{subeq:MF-4-7-0-6}), we clearly observe that two polynomial expressions, namely, Eqs.(\ref{subeq:MF-4-7-0-2}) and (\ref{subeq:MF-4-7-0-4}), contain no constant term compared to other expressions. This implies that there exist significant difference among the corresponding growth operations. Keep it in mind, let us recall operations introduced in Section 3. We indeed find that two types of operations associated with Eqs.(\ref{subeq:MF-4-7-0-2}) and (\ref{subeq:MF-4-7-0-2}), i.e., Type-I growth operation and V-fractal operation, are defined in a fashion independent of degree of vertex. On the contrary, the other four operations do greatly depend on degree of vertex. In a nutshell, we reach the following conjecture.

\textbf{Conjecture} \emph{ For a given tree $\mathcal{T}$ and a primitive growth operation $\mathcal{O}$, the formula for Wiener index of the resulting tree has no constant term when operation $\mathcal{O}$ is described regardless of degree of vertex. On the other hand, constant term is observed in the formula when operation $\mathcal{O}$ is described in a manner closely related to degree of vertex.}

\section{Applications}

In this section, we in depth discuss applications based on results above to some classic growth tree networks with interesting properties, such as, T-graph \cite{Agliari-2008}, Vicsek fractal \cite{Ma-EPL-2021} and Cayley tree \cite{Erzana-2020}. It should be mentioned that we study more general forms. In other words, those well-known trees are just specific examples of the following models. In particular, we obtained exact solutions of mean hitting time for random walks on these models. At the same time, mean shortest path length of two families of famous random tree networks, BA-scale-free tree and random uniform growth tree, are considered in detailed, and we then derive the analytic solutions. In addition, two variants associated with Wiener index on tree are also discussed in detail.

\subsection{Subdivision tree $\mathcal{T}^{m}(t)$}

As the first example model, subdivision tree $\mathcal{T}^{m}(t)$ is iteratively generated based on an arbitrary tree $\mathcal{T}$ by using $m$-order subdivision operation.

\textbf{Proof of Proposition I.7} This is proved by using Eqs.(\ref{eqa:MF-2-3-3}), (\ref{eqa:MF-3-1-1}) and (\ref{subeq:MF-4-7-0-1}). \qed

\subsection{Line graph}

First of all, we introduce the concept of line graph as follows. Given a graph $\mathcal{G}(\mathcal{V},\mathcal{E})$, one can obtain its corresponding \emph{line graph}, denoted by $\mathcal{G}_{\mathcal{L}}(\mathcal{V}_{\mathcal{L}},\mathcal{E}_{\mathcal{L}})$, that has as vertices the edges of $\mathcal{G}$, two edges being adjacent if they have an end in common \cite{Bondy-2008}.

\textbf{Proof of Proposition I.8} From definition of line graph, it is clear to the eye that each edge $e_{uv}$ in tree $\mathcal{T}$ corresponds to a unique vertex $w_{e_{uv}}$ in the line graph $\mathcal{T}_{\mathcal{L}}$. At first sight, this seems to be closely related to subdivision on edge mentioned in Section 3.1. Indeed, there exists a relationship that enables calculation of Wiener index $\mathcal{W}_{\mathcal{L}}$ of line graph $\mathcal{T}_{\mathcal{L}}$. Specifically, distance $d_{w_{e_{uv}}w_{e_{xy}}}$ between vertices $w_{e_{uv}}$ and $w_{e_{xy}}$ in line graph $\mathcal{T}_{\mathcal{L}}$ is equal to half of distance $d'_{w_{e_{uv}}w_{e_{xy}}}$ of the corresponding pair of vertices in set $\Lambda^{1}(1)$. The latter has been derived in Eq.(\ref{eqa:MF-4-1-2}) where parameter $m$ is assumed to be $1$, which implies

\begin{equation}\label{eqa:MF-5-2-1}
\mathcal{W}_{\mathcal{L}}=\frac{1}{2}\sum_{w_{e_{uv}},w_{e_{xy}}\in\mathcal{T}_{\mathcal{L}}}d_{w_{e_{uv}}w_{e_{xy}}}=\frac{1}{2}\mathcal{W}_{\mathcal{T}^{1}(1)}(2)=\mathcal{W}_{\mathcal{T}}-\left(
                                                             \begin{array}{c}
                                                               n \\
                                                               2 \\
                                                             \end{array}
                                                           \right).
\end{equation}
This completes the proof of Proposition I.8. Note that the same result as in Eq.(\ref{eqa:MF-5-2-1}) has been derived using another method \cite{Knor-2013}. \qed

\subsection{Classic BA-scale-free tree}

In 1999 \cite{Albert-1999-1}, Barab\'{a}si \emph{et al} revealed the scale-free feature popularly observed in a wide range of complex networks, and proposed the well-known BA-scale-free model $\mathcal{G}_{BA}(t)$ through two mechanisms, i.e., \emph{preferential attachment} and \emph{growth}. Roughly speaking, model $\mathcal{G}_{BA}(t)$ is iteratively built based on a seed $\mathcal{G}$ as follows: (i) at time $t$, a new vertex $u$ with $m$ edges is added into model $\mathcal{G}_{BA}(t-1)$, and (ii) each existing vertex $v$ in model $\mathcal{G}_{BA}(t-1)$ is connected to vertex $u$ with a probability $\prod_{v}$ proportional to its degree $k_{v}(t-1)$, namely, $\prod_{v}=\frac{k_{v}(t-1)}{\sum_{w\in\mathcal{G}_{BA}(t-1)}k_{w}(t-1)}$. Here, we are mainly interested in a special case of model $\mathcal{G}_{BA}(t)$ where parameter $m$ is assumed to equal $1$ and an edge serves as the seed. For convenience, the resulting model after $t$ time steps is denoted by $\mathcal{T}_{SF}(t)$ that is the so-called BA-scale-free tree. As well known, tree $\mathcal{T}_{SF}(t)$ obeys power-law degree distribution $P(k)\sim k^{-\gamma}$ where exponent $\gamma$ equals $3$. Below we provide an analytical formula of mean shortest path length $\langle\mathcal{W}_{SF}(t)\rangle$ of tree $\mathcal{T}_{SF}(t)$ by means of some consequences derived in Section 4.

\textbf{Proof of Proposition I.9} First, it is easy to see that there are $t+2$ vertices in BA-scale-free tree $\mathcal{T}_{SF}(t)$. Due to detailed description of tree  $\mathcal{T}_{SF}(t)$, Wiener index $\mathcal{W}_{SF}(t)$ is analyzed in an iterative manner as above. And then, mean shortest path length $\langle\mathcal{W}_{SF}(t)\rangle$ is obtained via Eq.(\ref{eqa:MF-2-2-2}). To this end, we need to build up a recurrence between quantities $\mathcal{W}_{SF}(t-1)$ and $\mathcal{W}_{SF}(t)$ of two consecutive growth trees $\mathcal{T}_{SF}(t-1)$ and $\mathcal{T}_{SF}(t)$, which is shown in the following form.

By definition in Eq.(\ref{eqa:MF-2-2-1}), the addition of new vertex $u$ into tree $\mathcal{T}_{SF}(t-1)$ has no influence on Wiener index $\mathcal{W}_{SF}(t-1)$. One can understand from construction of tree $\mathcal{W}_{SF}(t)$ that on average, each vertex $v$ with degree $k_{v}(t-1)$ in tree $\mathcal{W}_{SF}(t-1)$ is connected to vertex $u$ with probability $\prod_{v}=\frac{k_{v}(t-1)}{2(t+1)}$. We now need to determine distance $d'_{uv}(t)$ between vertex $u$ and each vertex $v$ in tree $\mathcal{T}_{SF}(t-1)$. As known, tree $\mathcal{T}_{SF}(t-1)$ is of stochastic form. Therefore, we get around this issue by making an expected estimation of quantity $\mathcal{W}_{SF}(t)$. In consequence, we can have

\begin{equation}\label{eqa:MF-5-3-1}
\begin{aligned}\mathcal{W}_{SF}(t)&=\mathcal{W}_{SF}(t-1)+t\sum_{v\in\mathcal{T}_{SF}(t-1)}\Pi_{v}+\sum_{v,w\in\mathcal{T}_{SF}(t-1)}\Pi_{v}\Pi_{u}+\sum_{v\in\mathcal{T}_{SF}(t-1)}\left(\Pi_{v}\right)^{2}\\
&\quad+\frac{1}{2}\sum_{v,w\in\mathcal{T}_{SF}(t-1)}(\Pi_{v}+\Pi_{w})d'_{vw}(t)+\frac{1}{2}\sum_{v,w\in\mathcal{T}_{SF}(t-1)}\Pi_{v}\Pi_{w}d'_{vw}(t)
\end{aligned},
\end{equation}
where $d'_{vw}(t)$ is distance between vertices $v$ and $w$ in tree $\mathcal{T}_{SF}(t)$. Note also that the last three terms in the first line are easy to derive, and two terms in the second line can be obtained in a similar way as in Eqs.(\ref{eqa:MF-4-5-4}) and (\ref{eqa:MF-4-5-10}), respectively. Hence, we omit the detailed proof here for the purpose of readability. Last but most importantly, we can see $\langle\mathcal{W}_{SF}(t)\rangle=O(\ln t)$ in the limit of large graph size, which is the same as some previous results derived using different methods \cite{Newman-2001,B-B-2004}. \qed

\subsection{Scale-free tree $\mathcal{T}^{\mathrm{II}}_{m}(t)$}

Inspired by the generative mechanisms behind BA-scale-free model $\mathcal{G}_{BA}(t)$, a great number of networked models with scale-free feature, both deterministic and stochastic, have been proposed and well studied in the past \cite{Newman-2010}-\cite{MA-W-W-L-pre-2020}. Among of them, scale-free trees $\mathcal{T}^{\mathrm{II}}_{m}(t)$ are recursively established through applying Type-II growth operation to an arbitrary tree $\mathcal{T}$. In particular, a single edge is always chose as a seed for convenience in the literature \cite{Chelminiak-2011}. Note that we consider more general models $\mathcal{T}^{\mathrm{II}}_{m}(t)$ below. After some simple arithmetic, trees $\mathcal{T}^{\mathrm{II}}_{m}(t)$ turn out to obey power-law degree distribution with exponent $\gamma=1+\frac{\ln(1+2m)}{\ln(1+m)}$.

\textbf{Proof of Proposition I.10}  This is proved by using Eqs.(\ref{eqa:MF-2-3-3}), (\ref{eqa:MF-3-5-1}) and (\ref{subeq:MF-4-7-0-5}). \qed

\subsection{Random uniform growth tree}

It is well known that random graphs, such as ER-model \cite{P-E-1959}, have attracted more attention in the last as they are believed to be potential candidates for modeling growth networks. While the majority of these models are proved to be relatively unreasonable with respect to some measures, for instance, degree distribution, it is of great interest to uncover some topological properties planted on them from the theoretical point of view. For example, average degree has a significant influence on phase transition of ER-model \cite{P-E-1959}. As such, we consider a class of random uniform growth trees $\mathcal{T}_{RG}(t)$ that are built in the following form. An arbitrary tree $\mathcal{T}$, denoted by $\mathcal{T}_{RG}(0)$, serves as seed. At each time $t$, a vertex $v$ is added into tree $\mathcal{T}_{RG}(t-1)$ and connected to each existing vertex $u$ with a probability $\prod$. As opposed to BA-scale-free tree $\mathcal{T}_{SF}(t)$, probability $\prod$ is equal to $1/|\mathcal{T}_{RG}(t-1)|$. It is easy to show that after $t$ times, the resulting tree $\mathcal{T}_{RG}(t)$ follows exponential degree distribution. In the following, we are mainly interested in mean shortest path length $\langle\mathcal{W}_{\mathcal{T}_{RG}}(t)\rangle$ on tree $\mathcal{T}_{RG}(t)$. Note also that a single edge is selected as the seed for generating tree $\mathcal{T}_{RG}(t)$.

\textbf{Proof of Proposition I.11} By analogy with the proof in Subsection 5.3, it is clear to understand that there are also $t+2$ vertices in random uniform growth tree $\mathcal{T}_{RG}(t)$ in total. With an in spirt similar manner, Wiener index $\mathcal{W}_{RG}(t)$ of growth tree $\mathcal{T}_{RG}(t)$ is derived. After that, mean shortest path length $\langle\mathcal{W}_{RG}(t)\rangle$ is given upon definition  Eq.(\ref{eqa:MF-2-2-2}). Toward this end, we first need to establish a connection of quantity $\mathcal{W}_{RG}(t)$ to $\mathcal{W}_{RG}(t-1)$. Due to concrete description of growth tree $\mathcal{T}_{RG}(t)$, this connection is expressed in the following form

\begin{equation}\label{eqa:MF-5-5-1}
\begin{aligned}\mathcal{W}_{RG}(t)&=\mathcal{W}_{RG}(t-1)+\frac{1}{t+1}\left[2\mathcal{W}_{RG}(t-1)+t+1+2\left(
                                                             \begin{array}{c}
                                                               t+1 \\
                                                               2 \\
                                                             \end{array}
                                                           \right)\right]\\
                                                           &\quad+\frac{1}{(t+1)^{2}}\left[\mathcal{W}_{RG}(t-1)+2\left(
                                                             \begin{array}{c}
                                                               t+1 \\
                                                               2 \\
                                                             \end{array}
                                                           \right)\right]\\
&=\left(1+\frac{1}{t+1}\right)^{2}\mathcal{W}_{RG}(t-1)+t+1+\frac{t}{t+1}
\end{aligned}.
\end{equation}
Based on the above equation, the correctness of Eq.(\ref{eqa:MF-5-5-0}) is consolidated in an iterative fashion. \qed

\subsection{Exponential tree $\mathcal{T}^{\mathrm{I}}_{m}(t)$}

Similarly, the corresponding deterministic version of tree $\mathcal{T}_{RG}(t)$ is often called uniform growth tree, denoted by $\mathcal{T}^{\mathrm{I}}_{m}(t)$. Specifically, tree $\mathcal{T}^{\mathrm{I}}_{m}(t)$ is created based on an arbitrary tree $\mathcal{T}$ by iteratively implementing Type-I growth operation. As mentioned above, some prior works pay attention on discussion on a specific case in which the seed is just a single edge. However, we concern more general version $\mathcal{T}^{\mathrm{I}}_{m}(t)$ whose seed is not necessarily an edge but an arbitrary tree.

\textbf{Proof of Proposition I.12} This is proved by using Eqs.(\ref{eqa:MF-2-3-3}), (\ref{eqa:MF-3-2-1}) and (\ref{subeq:MF-4-7-0-2}). \qed

\subsection{Generalized T-fractal $\mathcal{T}^{\mathrm{T}}_{m}(t)$}

In \cite{Redner-2001}, Redner \emph{et al} proposed the famous T-graph and have in depth discussed some structural parameters. In what follows, we will study the generalized T-fractal $\mathcal{T}^{\mathrm{T}}_{m}(t)$ that is iteratively generated upon an arbitrary tree $\mathcal{T}$ by performing T-fractal operation. It should be noted that the general model $\mathcal{T}^{\mathrm{T}}_{m}(t)$ is surely the well-known T-graph when we choose an edge as seed. Some variants associated with T-graph have used in wide range of applications, such as, Peano basin fractal \cite{Bartolo-2009}.

\textbf{Proof of Proposition I.13}  This is proved by using Eqs.(\ref{eqa:MF-2-3-3}), (\ref{eqa:MF-3-3-1}) and (\ref{subeq:MF-4-7-0-3}). \qed

\subsection{Generalized V-fractal $\mathcal{T}^{\mathrm{V}}_{m}(t)$}

As one of underlying models modeling regular hyperbranched polymers in chemistry, the classic V-fractal has been proposed by Vicsek in \cite{Vicsek-1983} and proven useful in some applications. More generally, we can obtain the generalized V-fractal $\mathcal{T}^{\mathrm{V}}_{m}(t)$ from an arbitrary tree $\mathcal{T}$ as seed by iteratively performing V-fractal operation. Particularly, if the seed is a single edge, then, after $t$ times, the resulting model is certainly the celebrated V-fractal. Here, we focus on random walks on generalized V-fractal $\mathcal{T}^{\mathrm{V}}_{m}(t)$. As a result, the previous result on V-fractal is easily obtained.

\textbf{Proof of Proposition I.14}  This is proved by using Eqs.(\ref{eqa:MF-2-3-3}), (\ref{eqa:MF-3-4-1}) and (\ref{subeq:MF-4-7-0-4}). \qed

\subsection{Generalized Cayley tree $\mathcal{T}^{\mathrm{III}}_{m}(t)$}

Another fundamental model, i.e., Cayley tree, has found a great variety of applications. For example, the classic Cayley tree whose seed a star is often used to model dendrimers in chemistry and biology \cite{Wu-2012}. In this subsection, we intend to study more general models that are called generalized Cayley tree $\mathcal{T}^{\mathrm{III}}_{m}(t)$. The model $\mathcal{T}^{\mathrm{III}}_{m}(t)$ is also constructed in an iterative way as above. That is to say, model $\mathcal{T}^{\mathrm{III}}_{m}(t)$ is obtained from the preceding model $\mathcal{T}^{\mathrm{III}}_{m}(t-1)$ by executing Type-III growth operation. It is worth noticing that the seed of model $\mathcal{T}^{\mathrm{III}}_{m}(t)$ is an arbitrary tree $\mathcal{T}$.

\textbf{Proof of Proposition I.15} This is proved by using Eqs.(\ref{eqa:MF-2-3-3}), (\ref{eqa:MF-3-6-1}) and (\ref{subeq:MF-4-7-0-6}). \qed

\subsection{Extended Wiener index based on multiplicative degree}

As stated in subsection 2.2, Wiener index $\mathcal{W}_{\mathcal{G}}$ of graph $\mathcal{G}$ is based on distance $d_{uv}$ between two vertices $u$ and $v$. Additionally, some variants have been defined in the rich literature \cite{Chen-2007}-\cite{Gutman-1996}. Below we focus on one of them, i.e., multiplicative degree Wiener index $\mathcal{W}^{\ast}_{\mathcal{G}}$. Specifically, quantity $\mathcal{W}^{\ast}_{\mathcal{G}}$ is the summation over multiplicative degree distances $d^{\ast}_{uv}$ of all possible pairs of vertices, say $u$ and $v$, in which $d^{\ast}_{uv}$ is defined as $k_{u}k_{v}d_{uv}$.

\textbf{Proof of Proposition I.16} This is an immediate consequence of Eq.(\ref{eqa:MF-4-5-10}), and we thus omit the detailed proof. Note that for a given graph $\mathcal{G}$, a more formula of parameter $\mathcal{W}^{\ast}_{\mathcal{G}}$ was reported in \cite{Chen-2007}. \qed

\subsection{Extended Wiener index based on additive degree}

As previously, another variant related to Wiener index $\mathcal{W}_{\mathcal{G}}$ of graph $\mathcal{G}$ is defined based on additive degree distance $d^{\dagger}_{uv}$ between vertex pair. That is to say, we denote by $\mathcal{W}^{\dagger}_{\mathcal{G}}$ the summation over additive degree distances $d^{\dagger}_{uv}$ of all possible pairs of vertices, say $u$ and $v$. Here, we study this quantity of a tree $\mathcal{T}$.

\textbf{Proof of Proposition I.17} This is an immediate consequence of Eq.(\ref{eqa:MF-4-5-4}), and we also omit the detailed proof.  Note also that for a given graph $\mathcal{G}$, a more formula of parameter $\mathcal{W}^{\dagger}_{\mathcal{G}}$ was reported in \cite{Gutman-1996}. \qed

\section{Some extremal problems and discussions}

As mentioned previously, tree as a special graph has been widely studied. In this study, we discuss some tree networks frequently observed in various kinds of fields. In particular, two fundamental and important structural parameters, Wiener index and mean hitting time, on these trees are studied in detail. It is well known that for a given tree $\mathcal{T}$ with $n$ vertices, the solution of Wiener index $\mathcal{W}_{\mathcal{T}}$  is subject to the following inequality

$$(n-1)^{2}\leq\mathcal{W}_{\mathcal{T}}\leq\left(
                                                             \begin{array}{c}
                                                               n+1 \\
                                                               3 \\
                                                             \end{array}
                                                           \right).$$
The equality on the left hand side of the above formula holds true when tree in question is a star. If tree under consideration is a path, then the equality on the right hand side of expression is achieved. Based on this, the mean shortest path length $\langle\mathcal{W}_{\mathcal{T}}\rangle$ of tree $\mathcal{T}$ follows

$$\frac{2(n-1)}{n}\leq\langle\mathcal{W}_{\mathcal{T}}\rangle\leq\frac{n+1}{3},$$
and, the mean hitting time $\langle\mathcal{H}_{\mathcal{T}}\rangle$ is immediately given by

$$\frac{2(n-1)^{2}}{n}\leq\langle\mathcal{H}_{\mathcal{T}}\rangle\leq\frac{(n+1)(n-1)}{3}.$$

Taken together, if a tree $\mathcal{T}$ has a more similar underlying structure to star, then the topological parameters, such as, Wiener index and mean hitting time, are closer to the corresponding lower bound. On the contrary, these parameters are closer to the associated upper bounds when the underlying structure of the tree is more like a path. This provides a guideline to create some iteratively growth tree networks that might be anticipated to show many other interesting properties, for instance, fractal feature. As discussed in the preceding sections, Vicsek fractal and T-graph are two representatives of such type of example trees partial because they possess fractal feature. Accordingly, they behave differently with both star and path in form. This leads to a fact that the solutions of their Wiener index and mean hitting time all fall into the scope bounded by lower and upper bounds. However, if we make a slight modification of growth way to produce Vicsek fractal and T-graph, the resulting trees will substantially deviate from the corresponding original tree and show more similar underlying structure to star. In other words, the exact solutions to Wiener index and mean hitting time approach the theoretical lower bound more rapidly. We also believe that one can have the ability to generate more trees whose Wiener index and mean hitting time are quite close to the lower bound using those constructive methods reported herein. More discussions are left for interested readers as an exercise. On the other hand, the problem of how to create tree models that not only show intriguing structural properties but also have Wiener index and mean hitting time much closer to the upper bound seems not easy to address. The intuition upon such an assertion is as follows. Using the analysis above, it is straightforward to see that trees of such type should possess almost the same underlying structure as path. In other words, the diameters of these trees should have the same magnitude of order as vertex number. As a result, it is impossible to observe some interesting features, which are prevalent in various kind of networked models, such as, fractal feature and scale-free feature, on trees of this kind. Note that, here, we have made use of some empirical analyses, such as (1) network $\mathcal{G}(\mathcal{V},\mathcal{E})$ that shows scale-free feature is often believed to have diameter $\mathcal{D}$ subject to $\mathcal{D}\sim \ln|\mathcal{V}|$ or even $\mathcal{D}\sim \ln\ln|\mathcal{V}|$ \cite{Cohen-2003} (In some extremal cases, one may also see scale-free networks following both $\mathcal{D}=\Theta(\ln|\mathcal{V}|)$ and $\mathcal{D}=o(|\mathcal{V}|)$ \cite{Ma-PRE-2021}), and (2) fractal networks $\mathcal{G}(\mathcal{V},\mathcal{E})$ always turn out to satisfy equalities $\ln\mathcal{D}=\Theta(\ln|\mathcal{V}|)$ \cite{Oliver-2013}.

\section{Conclusion}

To conclude, we study random walks on various kinds of recursive growth tree networks whose seed is not necessary a single edge but an arbitrary tree, which have found a wide range of applications in both theory and practice, and primarily consider mean hitting time. According to the elegant relationship between Wiener index and mean hitting time in Eq.(\ref{eqa:MF-2-3-3}), the exact solutions of mean hitting time on more general situation are derived in a series of combinatorial manners, which are based on the so-called Mapping Transformation, instead of using the commonly used spectral technique. The formulas obtained here completely cover the previously reported results in well-studied and simplest cases where an edge or a star is often used as a seed. From the theory point of view, this work establishes more general principle and enables ones to well understand underlying structure on recursive growth tree networks. Additionally, we extend the methods proposed to many other networked models and obtain analytical solutions of relevant parameters in a mathematically rigorous way. Last but most importantly, we also discuss some extremal problems in the realm of random walks on tree networks.

\section*{Acknowledgment}

The research was supported by the National Key Research and Development Plan under grant 2020YFB1805400 and the National Natural Science Foundation of China under grant No. 62072010.

\end{document}